\magnification=1200
\input amstex
\documentstyle{amsppt}
\pageheight{23truecm}
\pagewidth{16truecm}
\parskip 3pt plus 1pt minus 1pt
\hfuzz 0.5cm
\vbadness2000

\NoBlackBoxes


\def\const{\operatorname{const.}}

\def\A{{\Cal A}}
\def\B{{\Cal B}}
\def\C{{\Cal C}}
\def\D{{\Cal D}}
\def\E{{\Cal E}}

\def\G{{\Cal G}}

\def \O{{\Cal O}}

\def\S{{\Cal S}}
\def\T{{\Cal T}}
\def\U{{\Cal U}}

\def\IC{{\Bbb C}}
\def\IH{{\Bbb H}}

\def\IP{{\Bbb P}}

\def\IR{{\Bbb R}} 
\def\IS{{\Bbb S}}

\def\IZ{{\Bbb Z}}

\def\gG{{\frak G}}

\def \tq {{\tilde q}}
\def \tp {{\tilde p}}

\def \inv {P_{-1}}

\overfullrule=0pt

\def\ch{\raise 0.3ex\hbox{$\chi$}\kern-.15em}

\def\lto{\longrightarrow}



\def \BET {1}
\def \BETb {2}
\def \KET {3}
\def \KETb{4} 
\def \KLM {5}
\def \KLMb {6}
\def \KLR {7}
\def \KLS {8}
\def \KLT {9}
\def \NSU {10}
\def \SHS {11}
\def \SHSb {12}
\def \TAOa{13} 
\def \TAOb{14}
\def \TATa{15}
\def \TATb{16}
\def \UHL{17}


\topmatter
\title On the well-posedness of the wave map problem in high dimensions \endtitle
\leftheadtext{Nahmod, Stefanov and Uhlenbeck}
\rightheadtext {Wave Maps}
\author \ Andrea Nahmod, \ Atanas Stefanov \ and \ Karen Uhlenbeck \endauthor
\date    \enddate 
\address A. Nahmod, Department of Mathematics and Statistics, Lederle GRT,
University of Massachusetts, Amherst, MA 01003-4515
\endaddress
\email nahmod\@math.umass.edu \endemail
\address A. Stefanov, Department of Mathematics and Statistics, Lederle GRT,
University of Massachusetts, Amherst, MA 01003-4515
\endaddress
\email stefan\@math.umass.edu \endemail
\address K. Uhlenbeck, Department of Mathematics, The University of Texas at Austin,
Austin, TX 78712-1082\endaddress
\email uhlen\@math.utexas.edu\endemail
\thanks The first author was partially supported by NSF
grant DMS 9971159. This first and third authors acknowledge the support and hospitality of 
IAS at Princeton (May 2001), where part of this research was carried out. The third author would like to thank Montana State University for their hospitality as well.
\endthanks

\abstract We construct a gauge theoretic change of variables for the wave map from $\IR \times \IR^n$ into a compact group or Riemannian symmetric space, prove a new multiplication theorem for mixed Lebesgue-Besov spaces, and show the global well-posedness of a modified wave map equation - $n \ge 4$ - for small critical initial data. We obtain global existence and uniqueness for the Cauchy problem 
of wave maps into {\it compact} Lie groups and symmetric spaces with small critical initial data and $n \ge 4$. 

\endabstract

\endtopmatter

\document
\baselineskip=13pt

\head{0. Introduction }\endhead

The wave map equation between two Riemannian manifolds- the wave equation version
 of the evolution equations which are derived from the same geometric considerations as the harmonic map equation between two Riemannian manifolds-
has been studied by a number of mathematicians in the last decade.
The work of Klainerman and Machedon and Klainerman and Selberg \cite{\KLM} 
\cite{\KLMb} \cite{\KLS} studying the Cauchy problem for regular data 
is probably the best known. The more recent work of Tataru \cite{\TATa}, \cite{\TATb} and Tao \cite{\TAOa} \cite{\TAOb} relies 
and further develops deep ideas from harmonic analysis - in Tao's case 
in conjunction with gauge theoretic geometric methods - and thus 
seems very promising. Keel and Tao studied the one (spatial) dimensional 
case in \cite{\KETb}.

In \cite{\TAOa}, Tao established the global regularity  
for wave maps from $\IR \times \IR^n$ into the sphere $\IS^m$ when $n \ge 5$. Similar results to those of Tao were obtained by Klainerman and Rodniansky 
\cite{\KLR} for target manifolds that admit a bounded parallelizable structure.  

In this paper we are interested in revisiting this work. We study the Cauchy problem for wave maps from $\IR \times \IR^n$ into a (compact) Lie group (or 
Riemannian symmetric spaces) when $n \ge 4$ and establish global exisitence and uniqueness provided the Cauchy initial data are small in the critical norm. 
Similar results were obtained by Shatah and Struwe at roughly the same time when the 
target is any complete Riemannian manifold with bounded curvature. 

Our method combines both delicate  techniques from harmonic analysis with 
fairly standard global gauge theoretic geometric methods. Both our work and that of Shatah-Struwe \cite{\SHS} use the same gauge change; the analytic approach 
however, is significatly different as Shatah-Struwe base their results on Lorentz spaces and we use Besov spaces. Besov spaces are contained in Lorentz spaces -for appropriate indeces- (c.f. \cite{\SHSb} for example). Lorentz spaces seem to be more useful due to their better behaviour under coordinate transformations. 

It is interesting to note that in none of the works above is possible to obtain 
(strong) well posedness at the critical level for the wave map itself.
In other words even though one indeed has well posedness for the {\it gauged map}; 
there are no estimates available on {\it differences} for the {\it original wave map} itself and 
one cannot obtain any continuous dependence of the map on the data in the coordinate setting. It thus seems reasonable to think that the notion of wellposedness is not appropriate for this type of geometric equations at the critical level. The problem stems in that {\it well posedness} is not a gauge invariant notion; it is not even necessarily true that uniqueness in one coordinate system implies uniqueness in another directly. 

The plan for this paper is as follows. In section 1 we describe the geometry which translates wave maps into compact groups and Riemannian symmetric spaces to a gauge equation - the gauged wave map (GWM) -. This equation is overdetermined and we give a modified version (MWM). Section 2 containes the basic estimates for our theory, which involve multiplication theorems in Lebesgue-Besov spaces. Proposition 2.12 is the key estimate. This is the tool which contributes to handling the notoriously difficult first derivative non-linearity of the wave equation $$\square u + a \cdot du =0. $$ We obtain our results using the quadratic structure of the definition of $a$ in terms of $b$ ( which is linear in $du$ )
$$ \Delta a + div ( [ b,b ] + [a, a] ) =0. $$ This estimate is the subject of section 3. Section 4 contains the proof of the global well-posedness of the modified wave map equation (MWM) for small initial data in the scale invariant norm 
${\dot H}^{n/2}$. In section 4 we briefly outline the translation back to the original wave map coordinates. Our main result is the existence and uniqueness of global wave maps into compact Lie groups and symmetric spaces for small initial data in ${\dot H}^{n/2} \times {\dot H }^{n/2-1}$ for $n \ge 4$.

There are small difficulties in handling the case of non-compact symmetric spaces. The natural isometric embeddings are into spaces with indefinite metrics. For the standard methods on density theorems and coordinate changes to apply, it is necessary to know the existence of a Nash embedding into an Euclidean space with bounded geometry. 

Our results extend the results of Tao and Tataru for $M= \IS^m$. The Shatah-Struwe methods using Lorentz spaces are stronger since they obtain estimates for solutions with {\it variable} curvature. (On the surface, our difficulties with {\it non-compact} targets have somehow been circumvented in their work \cite{\SHS}).   

We have stated the results in sections 1 and 2 in great generality in the hopes that they may be applicable to other non-linear wave equations. The Appendix also contains alternate proofs of two of the multiplication theorems contained in the {\it Main Multiplication} theorem. These are the principal ones that are needed in 
estimating the non-linearity term in Theorem (2.13). 

\vskip .05in

The authors particularly thank both J. Shatah and T. Tao for their generous sharing of information and suggestions in a field relatively new to us, as well as for their enthusiastic support. 

\head{1. Formulation of the problem and gauge choices}\endhead

We regard the wave map equation as an equation given through covariant derivatives. These arise as follows: $$ s \,: \IR \times \IR^n \to M $$ where $M$ is an arbitrary Riemannian manifold and $$d \, s \, :  \T( \IR \times \IR^n ) \to \T M $$ where $\T(\IR \times \IR^n) = (\IR \times \IR^n ) \times ( \IR \oplus \IR^n)$.

\noindent Let $s^{\ast} \nabla $ be the pullback of the Levi-Civita connection on $M$ to  $s^{\ast} \T M $ via the map $s$. Then, in coordinate free notation, the wave map equation is

$$ s^{\ast} \nabla_0 \frac{\partial s}{\partial t} - \sum_{j=1}^n  s^{\ast} \nabla_j \frac{\partial s}{\partial x^j}\,=\,0.$$

Since the Levi-Civita connection on $M$ is torsion free, 

$$ s^{\ast} \nabla_j \frac{\partial s}{\partial x^k}=  s^{\ast} \nabla_k \frac{\partial s}{\partial x^j }\, , $$ for $j=0, 1, \dots, n \, , \, \, k=1, \dots, n$ where we have set $t = x^0$. 

We assume the map $s$ is topologically trivial which is usually implied by the later curvature bounds. (Note that the wave map fixes spatial infinity so topologically, $s: \IR^n \cup {\infty} \to M$). Hence,  $s^{\ast} \T M $ is the trivial bundle $ (\IR^{1} \times \IR^n ) \times \IR^m $. We also have control on the curvature of $ s^{\ast} \nabla $ via the equation 
$$[ s^{\ast} \nabla_j,  s^{\ast} \nabla_k ] = R(s) ( \frac{\partial s}{\partial x^j}, \frac{\partial s}{\partial x^k} ). $$

Our first theorem asserts that under smallness assumptions on $s \in L^{\infty}_t {\dot W}_x^{1, n/2}$, there is a unique choice of coordinates for $s^{\ast} \T M $. Given a smooth map $s$ with sufficient decay in asymptotics (to a point) at infinity, the initial coordinates  can be found by a partition of unity. The theorem we need is stated in a more general framework, as we hope to find applications for this theorem in gauge theory. 

\proclaim{ (1.1) Theorem } Let $ d + A$ be a smooth connection with compact structure groups $G$ over $\IR \times \IR^n$ or $I \times \IR^n$. Assume $A \sim 0$ at spatial infinity and let $F_A = d A + [ A, A] $ be the space-time curvature. Then there exists a positive constant $\epsilon = \epsilon( n, G) $ such that if the mixed space-time Lebesgue norm 
$$ \Vert F_A \Vert_{L^{\infty}_t L^{n/2}_x } < \epsilon \, , $$ then, there exits a unique smooth gauge change $g$, $g \sim I $ at spatial infinity, such that if 
$ {\tilde A} = g A g^{-1} - d g \, g^{-1} $ we have,  

\vskip .1in 

\roster
\item $\Vert  {\tilde A} \Vert_{{L_t^\infty} {\dot W}_x^{1, n/2}} \le c(n, G) \Vert F_A \Vert_{{L_t^\infty} {L_x}^{n/2} } $

\vskip .1in 

\item $ \sum_{j=1}^n \frac{\partial}{\partial x^j} {\tilde A}_j = 0 $  

\endroster

\endproclaim 

\demo{Proof} The method of proof follows the method used by the third author 
in \cite{\UHL}. We omit the dependence of constants on $G$ in the following proof. First, we fix each time slice $t=t_0$. The methods of  \cite{\UHL} show that in every ball $B_N = \{ x : |x| \le N \}$ there exists a gauge change $g_N$ such that the spatial part of the connection $ A_N = g_N A g_N^{-1} - d g_N \, g_N^{-1} $ satisfies 
$$ \Vert A_N \Vert_{t, {\dot W}_x^{1, n/2} (B_N) } \le \, c(n) \, \Vert F_A \Vert_{t, { L_x}^{ n/2} (B_N) }. $$ By taking $N \to \infty$,  
$\,  g_N  {\rightharpoonup} g $, $\, A_N \to {\tilde A}$, and we obtain a solution
$$\tilde A  = g A g^{-1} - dg \, g^{-1}, \qquad \sum_{j=1}^n \frac{\partial}{\partial x^j} {\tilde A} = 0 $$ on each time slice $(t, \IR^n)$ which satisfies on all $ \IR^n$, $$\Vert {\tilde A} \Vert_{t, {\dot W}_x^{1, n/2} }\le c(n) \Vert F_A \Vert_{t, L_x^{n/2}}\,\,  . $$ Since $A$ is asymptotic to $0$ at infinity, $d g$  is as well, and we may choose $g \sim I$ at spatial infinity as well. Let $g = \text{exp}( u ) $. We fix a time slice and then differentiate in $t$. Namely, if 
$$ \sum_{j=1}^n \frac{\partial}{\partial x^j} \text{exp}(u) {\tilde A}_j  \text{exp}(-u) - \frac{\partial}{\partial x^j}  \text{exp}(u)   \text{exp}(-u) =0 $$ is the equation at the time slice $t_0$, the derivative at $t=t_0$ is 
$$ L_{\tilde A} u = \Delta u + d\ast [u, {\tilde A} ] = \Delta u + [ du, {\tilde A}].$$ Here we use the fact that $\dfrac{\partial}{\partial x^j} {\tilde A}_j=0$ at $t=t_0$.  Examine the properties of this linear map -which is the derivative- $$ L_{\tilde A} : {\dot L}^{ 2, n/2} \to L^{n/2}, \qquad L_{\tilde A} = \Delta + \text{lower order terms}.$$ We have $$\align \Vert [ du, {\tilde A}]\Vert_{L^{n/2}} &\le 2 \Vert du \Vert_{L^n} \Vert{\tilde A} \Vert_{L^n} \\
&\le c(n) \Vert u \Vert_{{\dot W}^{2,n/2}} \Vert A \Vert_{{\dot W}^{1, n/2}}\\
&\le c(n) \, \epsilon  \, \Vert u \Vert_{{\dot W}^{2,n/2}}. \endalign $$ Choose $\epsilon$ so that $ c(n) \, \epsilon < 1/2$; we have that the lower order term is small enough for $L_{\tilde A}$ to be invertible.  The precise estimate is for 
$ u = \dfrac{ \partial g }{\partial t }  g^{-1} $ where 
$$L_{\tilde A} u + \sum_{j=1}^n \frac{ \partial }{\partial x^j } ( g \, \frac{\partial A}{\partial t} \, g^{-1} ) = 0 .$$ Note that an estimate on $\frac{\partial g}{\partial t}$ is available by the general methods we have been using. 

To obtain an estimate on the time component ${\tilde A}_0 = g A_0 g^{-1} - \dfrac{\partial g}{\partial t} g^{-1} $, note that 
$$ d {\tilde A}_0 - \frac{\partial {\tilde A}}{\partial t} + [ \tilde A, A_0 ]= (F_A)_{(\text{space,time})}. $$ Since $d\ast {\tilde A} =0 $, we have $$ \Delta {\tilde A}_0 + \frac{\partial}{\partial x^j}[{\tilde A}_j, A_0] = \sum_{j=1}^n  \frac{\partial}{\partial x^j} ( F_{\tilde A})_{j,0}.$$ Let $\mu_j = {\Delta}^{-1/2} \frac{\partial}{\partial x^j}$. Then 
$$\align \Vert {\tilde A}_0 \Vert_{{\dot W}^{1, n/2}} &= \Vert \Delta^{1/2} {\tilde A}_0 \Vert_{{L}^{ n/2}} \\
&\le \Vert \sum_{j=1}^n \, \mu_j \, [{\tilde A}_j, {\tilde A}_0]_{L^{n/2}} + \Vert \sum_{j=1}^n \, \mu_j (F_{\tilde A})_{j,0} \Vert_{L^{n/2}} \\
&\le c(n)\, ( \Vert {\tilde A} \Vert_{L^n} \Vert {\tilde A}_0 \Vert_{L^n} + \Vert F_{\tilde A} \Vert_{L^{n/2}} )\\
&\le c(n) \, ( \, \epsilon \, \Vert {\tilde A}_0 \Vert_{L^n} + \Vert F_{A} \Vert_{L^{n/2}} ) .\endalign $$ Again, under the assumption $c(n) \, \epsilon < 1/2 $, we have
$$ \Vert {\tilde A}_0 \Vert_{{\dot W}^{1, n/2}} \le 2 c(n) \,  \Vert F_{A} \Vert_{L^{n/2}}  $$ as claimed. \qed

\enddemo

\proclaim{ (1.2) Corollary } Theorem (1.1) remains true if $A \in L^{\infty}_t {\dot W}^{1, n/2}_x$ and $F_A \in  L^{\infty}_t {L}^{n/2}_x$ 
\endproclaim 

\demo{Proof} Approximate $A$ by smooth connections $A_{\alpha} \to A$ in ${\dot W}^{1, n/2}$. For each $\alpha$, construct $g_{\alpha}$ as in the theorem and $A{\alpha} = g_{\alpha} A_{\alpha} g_{\alpha}^{-1} - {\hat d} g_{\alpha} g_{\alpha}^{-1}$ a space-time connection $1$-form which satisfies the estimates. We have denoted by $\hat d= ( \dfrac{\partial}{\partial t}, d) $ the full derivative. Then, 
$$\align  \Vert {\hat d} g_{\alpha} \Vert_{L^n} &=  \Vert {\hat d} g_{\alpha} \,  g_{\alpha}^{-1} \Vert_{L^n} \\
& \le  \Vert {\tilde A}_{\alpha} \Vert_{L^n} +\Vert g_{\alpha} {A}_{\alpha} g_{\alpha}^{-1} \Vert_{L^n} \\
& \le \Vert {\tilde A}_{\alpha} \Vert_{L^n} +\Vert  {A}_{\alpha}  \Vert_{L^n}. \endalign $$ Hence $ \Vert {\hat d} g \Vert_{L^n}$ is bounded on each time-slice. To complete the estimate note that 
$$\align (\frac{\partial g}{\partial x^j})_{\alpha} &= ({\tilde A}_{\alpha})_j g_{\alpha} - g_{\alpha} (A_{\alpha})_j \\
d\, (\frac{\partial g}{\partial x^j})_{\alpha} &= d\, ({\tilde A}_{\alpha})_j g_{\alpha} - g_{\alpha}\, d\, (A_{\alpha})_j +   ({\tilde A}_{\alpha})_j \, d\, g_{\alpha} - \, d\, g_{\alpha} (A_{\alpha})_j .\endalign $$ Then again, on each time slice $$ \Vert \hat{d} g_{\alpha} \Vert_{W^{1, n/2}} \le \Vert {\tilde A}_{\alpha} \Vert_{W^{1, n/2}} + \Vert {A}_{\alpha} \Vert_{W^{1, n/2}} + ( \Vert {\tilde A}_{\alpha} \Vert_{L^{n}} + \Vert {A}_{\alpha} \Vert_{L^{n}}) \, \Vert d g_{\alpha} \Vert_{L^n}  $$ is also bounded. 

In each time slice, we have subsequences which converge to weak limits 
$\, g_{\alpha'} \rightharpoonup g $ in $\, {\dot W}^{2, n/2}$. However, the weak limit is unique. Suppose not. Then,
$$  g_{\alpha'} \rightharpoonup g \,  ,  \qquad \qquad   g_{\alpha''}
\, \rightharpoonup \,  h g \, ; $$
$$ {\tilde A} = g A g^{-1} - d g \, g^{-1} \qquad \text{ and }  \qquad  {\tilde {\tilde A}}= h {\tilde A} h^{-1} - d h \, h^{-1}. $$ Both ${\tilde A}$ and ${\tilde {\tilde A}}$ satisfy the time-slice estimate 
$$ \Vert {\tilde A} \Vert_{{\dot W}^{1, n/2}} + \Vert {\tilde {\tilde A}} \Vert_{{\dot W}^{1, n/2}} \le \, 2 \, \epsilon $$ as well as
$$ \sum_{j=1}^n \frac{\partial}{\partial x^j} {\tilde A}_j = \sum_{j=1}^n \frac{\partial}{\partial x^j} {\tilde {\tilde A}}_j =0 .$$ But $ d h = h {\tilde A} - {\tilde{\tilde A}} \, h $ and $\Delta h = (d h {\tilde A}) - ({\tilde{\tilde A}} \, d h) $. If we let $k_n$ be the appropriate Sobolev constant, 
$$\align \Vert d h \Vert_{L^n} &\le \, c(n) \, \Vert \Delta h \Vert_{L^{n/2}} \\
& \le c(n) \, k_n\, \Vert d h \Vert_{L^n} \, (  \Vert {\tilde A} \Vert_{L^{n}} + \Vert {\tilde {\tilde A}} \Vert_{L^{n}}) \\
& \le c(n) \, k_n\, \epsilon \,  \Vert d h \Vert_{L^n}. \endalign $$ If  $2 \, c(n) \, k_n \, \epsilon < 1$ we have that $d h =0$. Since $ h \sim I$ at infinity, $ h \equiv I$. Thus the weak limit is unique. Hence, $g$ is unique and 
$$ {\hat d} g_{\alpha} \rightharpoonup {\hat d} g \, \text{ { in }}
{\dot W}^{1, n/2}\, , \qquad {\tilde A}_{\alpha} \rightharpoonup
{\tilde A} \, \text{ { in }} {\dot W}^{1, n/2}$$ and 
$$ \Vert {\tilde A} \Vert_{L^{\infty}_t{\dot W}^{1, n/2}_x} \le \, c(n) \, \Vert F_A \Vert_{L^{\infty}_t{L}^{n/2}_x}  $$ as claimed. \qed 
\enddemo

\proclaim{ (1.3) Corollary} Let $s: \IR^1 \times \IR^n \to M$  be an arbitrary map. Suppose the curvature $R(M)$ is bounded by $K$. There exists $0 < \delta = \delta( n, K) $ such that if $ s: \IR^n \cup \{ \infty \} \to N$ is topologically trivial, and $$ \Vert \, d s \Vert_{L^{\infty}_t L^n_x} \le \delta\, ,  $$ then there exists a unique frame in $s^{\ast} \T M $ such that the hypotheses of (1.1) - the main gauge-fixing theorem-  are satisfied. 

\endproclaim 

\demo{Proof} Since $$( F_{ s^{\ast} \nabla} )_{kj} = R( s ) ( \frac{\partial s}{\partial x^j}, \frac{\partial s}{\partial x^k} )\, , $$ on time slices we have the estimate $$\Vert F_{ s^{\ast} \nabla} \Vert_{L^{n/2}} \le K \Vert  d s \Vert^2_{L^n} .$$ The desired conclusion follows by choosing $ \delta >0$ such that $K \delta < \epsilon$ where $ 0 < \epsilon = \epsilon(n, M) $ is as in Theorem (1.1). \qed
\enddemo

\vskip .1in 
Next we give a coordinate invariant description of the wave equations. Let 

$$ D = s^{\ast} \nabla = d + a\, ,  $$ where the curvature of $d$ is $$ F_A = (R \circ s ) ( d\, s\, , \, d \, s ) . $$ The term   $(R \circ s )$ is not explicit unless one is working on a Lie group or symmetric space. 

Let $b = ds $. Then the equations themselves are written 
$$ D_0 b_0 - \sum_{j=1}^n D_j b_j = 0. $$ Because the Levi-Civita connection on $M$ has no torsion, we find 
$$ D_k b_j = D_j b_k \, , \quad k=0, 1, \dots, n \, , \, \, j=1, 2, \dots, n. $$ This is a non-linear first order hyperbolic system. It may be that the correct method is to analyze this directly. In keeping with the present standard methods, we convert it to a single equation using Hodge theory.

\proclaim{ (1.4) Theorem } Let $b = d \phi + d^{\ast} \psi$. Then the wave map equations can be rewritten as 

$$\align {\text{(a)}}\quad  &\square \phi  + ( a, b ) = 0 \\
{\text{(b)}}  \quad  &\square \psi  +  a \wedge  b  = 0 \\
{\text{(c)}}\quad  &b= d \phi + d^{\ast} \psi \\
{\text{(d)}} \quad &d a + [a, a] = R(x) [ b, b] \\
{\text{(e)}}\quad &\sum_{j=1}^n \frac{\partial}{\partial x^j} \, a_j = 0.  \endalign $$ Here $R(x)$ is the Riemannian curvature of $M$ evaluated at $s(x)$. 

\noindent The initial data on $\phi$ and $\psi$ can be taken to be 

$$\align &\phi(0, x) =0\, ,  \qquad \qquad \psi(0,x)=0 \\
& \frac{\partial \phi}{\partial t} (0, x) = b_0 (0, x) \qquad \qquad  \frac{\partial \psi}{\partial t}_{0,j} (0, x) = b_j (0, x)\\
& \frac{\partial}{\partial t} \psi_{j,k} (0, x) = 0 \qquad j,\,  k \neq 0  
\endalign $$

\endproclaim 

\demo{Proof} Let $b = \square q $ with $q(0,x) = \frac{\partial q}{\partial t}(0, x) =0 $ where $\square = d  \, d \ast +  d \ast \, d $ and $ d \ast= div_{\text(space, time)} $ computed using the Lorentz metric. 

Let $$\phi = d \ast q = \frac{\partial}{\partial t} b_0 - \sum_{j=1}^n \frac{\partial}{\partial x^j} b_j $$ and $$\psi = d q =\frac{\partial}{\partial t} b_j - \frac{\partial}{\partial x^j} b_0 \, , \,\frac{\partial}{\partial x^k} b_j - \frac{\partial}{\partial x^j} b_k. $$ Hence $$\square \phi = d \ast b \quad \text{ and } \quad \square \psi = d \wedge b $$ So, $ b = d \phi + div_{\text(space, time)} \psi $. Note $d \psi =0$ automatically. 

The initial data clearly consists of $\phi(0, x) =0 $, $\psi(0,x) = 0$. Hence, 
$$b_0 = \frac{\partial}{\partial t} \phi - \sum_{j=1}^n \frac{\partial}{\partial x^j} \psi_{j,0}\, ,  \quad {\text{ so } } \quad b_0(0,x) = \frac{\partial}{\partial t} \phi(0, x).  $$ Likewise,  
$$ b_j =  \frac{\partial}{\partial x^j} \phi - \frac{\partial}{\partial t} \psi_{0, j} + \sum_{k=1}^{n} \frac{\partial}{\partial x^k} \psi_{k, j }\, , \quad {\text{ so } } \quad b_j (0, x) = \frac{\partial}{\partial t} \psi_{0, j} (0, x). $$ Note also that 
$$\frac{\partial }{\partial t} \psi_{j,k} = \frac{\partial }{\partial x^j} \psi_{0,k} + \frac{\partial}{\partial x^k} \psi_{j,0}\, ,  \quad \text{ so } \quad \frac{\partial}{ \partial t} \psi_{j,k} (0, x) = 0. $$

\enddemo 
 The last equation (d) of (1.4) is not determined by the rest of the data since the curvature depends on the original map (and gauge change). No general formula is available. This would not preclude {\it a priori} estimates. However, the estimates 
for our global existence and uniqueness theorem for wave maps are done in Besov spaces (which here prove inferior to the Lorentz spaces). The equation 
$$ d a + [a, a] = R(x) [b, b] $$ however behaves `badly' (for bounded
$R(x)$) in this context. Hence we must restrict the manifold $M$ to a  group or a Riemannian symmetric space. 

\proclaim{ (1.5) Theorem} If $M = G$ or $M = H/G$ where $G$ is a compact Lie group, then the equation (d) in Theorem (1.4) can be replaced by the equation 
$$\text{ (d)' } \quad  d a + [a,a] + [b, b] = 0 .$$ Moreover, the original map $s: \IR^1 \times \IR^n \to G $ (or $H/G$ ) can easily be reconstructed from the fact that $\, d\, +\, a\, +\, b  \, \text{ and } \, d \, + \, a \, - \, b $ are flat connections. Let $$ ( d\, +\, a\, +\, b ) g^{+} = 0 \qquad \text{ and } \qquad  ( d\, +\, a\, -\, b ) g^{-} =0 . $$ The original map is $g = g^{+} \cdot g^{-}$. 
\endproclaim

\demo{Proof}
The computations for a Lie group are straightforward if we remember that $\T^{\ast} G = \gG$, the Lie algebra of G, that $[\cdot, \cdot]$ generates curvature, and that the structure group is a specialization of the orthogonal group. The symmetric space case is best understood by regarding $M$ as an $Ad$ orbit in the possibly non-compact group $H$. That is, 
$$ M = Ad H ( {\hat i}) $$ and $G$ is the (compact) isotropy subgroup of $\hat i$. 
For $\IH^m$, $H$ is the Lorentz group $\O(1,m)$ and $G$ is the Euclidean group $\O(m)$. Choose $ {\hat i} = {\text{diag}}( 1, -1, -1, \dots, -1)$. Then by construction $b$ will always lie in the off-diagonal vectors 
$$ b_j = \left[ \matrix 0 & \cdot & v_j & \cdot \\ \cdot & 0 & \dots & 0 \\
- v_j^{\ast} & 0 & \dots & 0 \\
\cdot & 0 & \dots & 0 \endmatrix \right]\, , $$ and the compact structure group $\O(m) $ is represented on the diagonal.  The construction cannot work for non-compact Lie groups $H$ such as $\O(1,m)$ since the do not have bi-invariant Riemannian metrics. \qed

\enddemo

\proclaim{ (1.6) Corollary} Suppose $M= G$ or $M= H/G$. Then a subset of the gauged wave map equations (a)--(e) (GWM) has a structure of a non-linear wave system of integral differential operators.
$$\align {\text{(a)}}\quad  &\square \phi  + ( a, b ) = 0 \\
{\text{(b)}}  \quad  &\square \psi  +  a \wedge  b  = 0 \\
{\text{(c)}}\quad  &b= d \phi + d^{\ast} \psi \\
{\text{(d)}} \quad &\Delta a_j + \sum_{j=1}^n \frac{\partial}{\partial x^k} [a_k, a_j] +  \frac{\partial}{\partial x^k} [b_k, b_j]\, , \quad j=0, 1, \dots n. \endalign $$

\endproclaim 

\demo{Proof} The wave equation structure of the system in $\phi$ ad
$\psi$ is clear, and $b$ is a (linear) first order derivative of
$\phi$ and $\psi$ (note that the initial data for $\phi$ and $\psi$
have been worked out as coupled to that of $b$.) $$\align \quad  &\phi(0, x) = \psi( 0, x) = 0  \\
 \quad  &\frac{\partial \phi}{\partial t} (0, x) = b_0 (0, x) \\
\quad  &\frac{\partial {\psi}_{j,0}}{\partial t} (0, x) = b_j (0, x) \\
\quad &\frac{\partial {\psi}_{j,k}}{\partial t} (0, x) = 0 .\endalign $$ To obtain the last equation, note that 
$$ \frac{\partial}{\partial x^k} a_j - \frac{\partial}{\partial x^j} a_k + [ a_j, a_k] + [b_j, b_k] = 0 $$ for $k=1, 2, \dots, n \, , \quad j=0, 1, \dots, n$ and $\frac{\partial}{\partial x^0} = \frac{\partial}{\partial t}$. Since $\sum_{k=1}^n \frac{\partial}{\partial x^k} a_k =0$ , we obtain our new equation by taking the divergence. We call this system the {\it modified wave map} (MWM). By indirect arguments, it is clear that if $( \phi, \psi, b, a)$ have initial data which satisfy appropriate constraints, the evolution, at least in the smooth case, will actually be data coming >from a wave map. The direct argument is not available to us however.

\enddemo

\head{2. Definitions and Product Estimates }\endhead

In this section we set out definitions, notations and basic estimates 
that will be used throughout the paper. We shall frequently use the notation 
$A \lesssim B$ to mean $ A \le \const B$ for some positive constant $\const $ which is allowed to vary from line to line but does not depend on any of the relevant parameters in the estimates.

We begin reviewing some Littlewood-Paley theory. Let $\phi(t, x)$ be a function on $\IR \times \IR^n$, we define the {\it spatial} Fourier transform $\hat{\phi}(t,
\xi)$ by 

$$\hat{\phi}(t, x) = \int_{\IR^n} e^{-2 \pi i x \cdot \xi} \phi( t, x) \, dx $$

We define now the usual Littlewood-Paley projection operators $P_k$ and $Q_k$. To that effect, let $m(\xi)$ be a non-negative radial bump function supported on the ball $|\xi |\leq 2$ and equal to $1$ on the ball $|\xi| \leq 1$. Then for each integer $k$ we define $P_k(\phi)$ the projection onto the frequency ball $|\xi| 
\lesssim 2^k \, $ by 
$$ \widehat{ P_k (\phi) }(\xi) := m(2^{-k} \xi ) \hat{\phi}(t, \xi ).  $$ 
Note that $ P_k \to 0$ in $L^2$ as $ k \to -\infty$ while $P_k \to I$ in $L^2$ as $k \to \infty$. 
\vskip .1in

\noindent The operator $Q_k$ is the projection  
onto the frequency annulus $|\xi| \sim 2^k \,$ given by the formula,
$$ Q_k := P_k  -  P_{k- 1}. $$
\noindent We note that if we let $ \psi({\xi}) := m( \xi ) - m (2 \xi ) $, 
then $\psi$ is supported on the annulus $ \, 1/2 \le |\xi | \le 2 \, $,  for all  
$\xi \neq 0, \quad \sum_{k \in \IZ} \psi(2^{-k} \xi) \equiv 1$, and 
$$ \widehat{ Q_k (\phi)}(t, \xi) = \psi( 2^{-k} \xi ) \hat{\phi}(t, \xi). $$
The Littlewood-Paley projections are bounded operators in all the Lebesgue spaces 
and commute with any constant coefficient differential operator. Finally 
we note that $Q_k$ is given by a convolution kernel whose $L^p$-norm equals $2^{ ( k n )( 1 - 1/p)}$ for all $ 1 \ge p \le \infty$. In particular its $L^1$- norm
is identically $1$ for all $ k \in \IZ$. 

Let $ j = 0 $ or $j = 1$ and let $k \in \IZ$.  Following \cite{\TAOa} and also \cite{\KLR}, we introduce $\S^{(-j)}_k (\IR \times \IR^n)$, the {\it Strichartz space at frequency } $2^k$ to be the space of functions whose space-time norm is given by: 
$$ \| \Phi \|_{\S^{(-j)}_{k}} := \sup_{ q, r \in \A} 2^{k (\frac{1}{q} + \frac{n}{r} - j) }\bigl ( \| \Phi\|_{L^q_t L^r_x} \, + \, 2^{-k}  \| \partial_t \Phi\|_{L^q_t L^r_x}\bigr), $$
where $\A := \{ (q,r) \, : \, 2 \le q,\, r \, \le \infty, \, \dfrac{1}{q} + \dfrac{n-1}{2r} \le \dfrac{n-1}{4} \, \} $ is the set of admissible Strichartz exponents. We remark that when $j=0$ the spaces above are $\dot{H}^{n/2}$- normalized and correspond to Tao's spaces $S_k$ in \cite{\TAOa}. We also note that for each $n \ge 4$,  only specific values of $(q, r)$ are needed. Finally observe that 
control of the $\S^{(-j)}_k$ norm gives, for example, the estimates:

$$  \| Q_k(\phi) \|_{L^2_t L_x^{\frac{2(n-1)}{(n-3)}}} + 2^{-k} \| \partial_t Q_k(\phi) \|_{L^2_t L_x^{\frac{2(n-1)}{(n-3)}}} \, \le \, 2^{k ( j + \frac{n}{(n-1)} - \frac{(n+1) }{2} ) } \, \| Q_k(\phi)\|_{\S^{(-j)}_k}
\tag2.1 $$

$$\| Q_k(\phi) \|_{L_t^{\infty} L_x^{2}} + 2^{-k} \| \partial_t Q_k(\phi) \|_{L^{\infty}_t L_x^{2}} \, \le \, 2^{k (j - \frac{n}{2}) } \, \| Q_k(\phi)\|_{\S^{(-j)}_k}\tag2.2 $$

$$  \| Q_k(\phi) \|_{L^2_t L_x^{\infty}} + 2^{-k} \| \partial_t Q_k(\phi) \|_{L^2_t L_x^{\infty}}  \, \le \, 2^{k ( j -  \frac{1}{2}) } \, \| Q_k(\phi)\|_{\S^{(-j)}_k}\tag2.3 $$

\noindent Finally we state the Strichartz estimates in this framework (c.f. \cite{\TAOa} \cite{\KET} and references therein ). 

\vskip .1in 

\proclaim{(2.4) Theorem (Strichartz Estimates) }

Let $k$ be an integer and let $\Phi$ be any function on $\IR \times \IR^n$ with spatial Fourier support on the annulus $|\xi| \sim 2^k$. Then 
$$ \| \Phi \|_{S_k^{(-j)}} \lesssim \| \Phi(0, \cdot) \|_{ \dot{H}_x^{n/2-j}} + \| \partial_t \Phi(0, \cdot) \|_{ \dot{H}_x^{n/2 -(j +1) }}+ 2^{k( \frac{n}{2} - ( j +1) )} \, 
\| \square \Phi \|_{L^1_t L_x^2}. $$  

\endproclaim 

\vskip .1in

\proclaim{ (2.5) Definition }

Let $\S^{(-j)}$ be the space of functions on $\IR \times \IR^n$ whose norm is given by 
$$ \| \phi \|_{\S^{(-j)}} := \bigl( \sum_{k \in \IZ} \| Q_k (\phi) \|_{S_k^{(-j)}}^2 \bigr)^{1/2}. $$   
\endproclaim

\proclaim{ (2.6) Definition} A pair $(q,r)$ is said to be  {\rm sharp 
admissible } if $2 \le q, r \le \infty$ and  
$$\frac{1}{q} + \frac{n-1}{2r} \, = \, \frac{n-1}{4}. $$ 
\endproclaim

\proclaim{Remark} If $n \ge 4$ and $(q, r)$
is sharp admissible then $s = 1/q + n/r - 1 > 0$. Also, in particular, 
$q \ge 2$ and $2 \le r \le \frac{ 2 (n-1)}{n-3}$.  

\endproclaim
 
\proclaim{ (2.7) Lemma} For any $ j \ge 0$ we have that 
$$ \sup_{ (q,r) - \text{admissible}} 2^{k ( \frac{1}{q}
+\frac{n}{r}-j)} \| Q_k(f) \|_{ L^q_t L^r_x} \, = \,  \sup_{ (q,r) -
\text{sharp admissible}} 2^{k ( \frac{1}{q}
+\frac{n}{r}-j)} \| Q_k(f) \|_{ L^q_t L^r_x}. $$
In other words, 
$$ \align & \| f \|_{\S^{(-j)}} \, = \\
&= \bigl( \sum_{k \in \IZ} \quad \bigl|\sup_{ (q,r) -
\text{sharp admissible}} 2^{k ( \frac{1}{q}
+\frac{n}{r}-j)} \bigl(  \| Q_k(f) \|_{ L^q_t L^r_x }+  2^{-k} \| \partial_t Q_k(f) \|_{ L^q_t L^r_x }\bigr) \, \bigr|^2\bigr)^{1/2}. \endalign $$ 
\endproclaim

\demo{Proof} Let $(q,r)$ be admissible but not sharp admissible and
define $r_0$ such that $(q, r_0)$ is {\it sharp} admissible. Then it
is clear that $r_0 < r$. Let $s > 0$ be such that $$\frac{1}{r} =
\frac{1}{r_0} - \frac{\gamma}{n}.$$ By the Sobolev embedding we then have
$$\| Q_k(f)\|_{L^q_t L^r_x} \le 2^{k \gamma} \| Q_k(f)\|_{L^q_t
L^{r_0}_x}.$$ 

Note that when $r= \infty$, $W^{n/{\gamma}, r_0} \hookrightarrow BMO$ but 
$\|Q_k(f)\|_{L^{\infty}} \sim \|Q_k(f)\|_{BMO}$. Hence {\it on each
Littlewood-Paley piece} the Sobolev embedding above holds.  

Then,
 $$2^{k( \frac{1}{q} +\frac{n}{r}-j)} \| Q_k(f) \|_{ L^q_t L^r_x} \leq
2^{k( \frac{1}{q} +\frac{n}{r}+ \gamma -j)} \| Q_k(f) \|_{ L^q_t L^{r_0}_x } \le  2^{k( \frac{1}{q}
+\frac{n}{r_0} -j)} \| Q_k(f) \|_{ L^q_t L^{r_0}_x} $$
and similarly for $2^{-k} \| \partial_t  Q_k(f) \|_{ L^q_t L^r_x} $ 
>from where the conclusion follows \qed. 

\enddemo

In what follows we will denote by  $P_{-1}: = \nabla \Delta^{-1}$ be the pseudodifferential operator defined by 
$$\widehat { P_{-1} f}(t, \xi) = \frac{1}{|\xi|} \hat{f}(t, \xi). $$

\proclaim{ (2.8) Definition}
We denote by $\B_p$ be the Banach space of functions on $\IR \times \IR^n$ 
 whose norm is given by 
$$ \| f\|_{\B_p} \, := \, \bigl( \sum_{k \in \IZ} \| Q_k (f) \|^p_{ L^1_t L^{\infty}_x} \bigr)^{1/p}$$  for $1 \le p < \infty$ and suitable modified with the $\ell^\infty$-norm when $p= \infty$

\endproclaim
\proclaim { Remark} Note that it follows naturally from the embeddings $\ell^p \subset \ell^q $ that $ \B_p \subset \B_q $ for $1 \le p < q \le \infty. $
\endproclaim 

\vskip .1in

We proceed to prove the {\it Main Multiplication
Estimate}. The point of it is that it implies in particular the 
three multiplication estimates that will be needed later and more. It thus
gives a unified framework under which to understand the action of the
`inverse gradient' $\inv$ on the space $S^{(-1)} \times
S^{(-1)}$. For  solutions of the homogeneous 
wave equation, Klainerman and Tataru \cite{\KLT} obtained 
the first bilinear estimates of this type on an improved range; 
those can be viewed as generalizations of the 
well-known Strichartz-Pecher inequalities. 

We first need some definitions. In what follows, for any $a \in \IR$, 
we will denote by $a^{-}$ and  $a^{+}$ the real number $a -
1/100$ and $ a + 1/100$ respectively. The constant
$1/100$ is of course arbitrary; any (fixed) small positive number will
do.  

\vskip .1in 

\proclaim{(2.9) Definition} Let us denote by $\C$, $\D$, $\E$,  $\G$ the
following sets of pairs $(q, p)$ where $q, p \ge 1$. 

$$\C = : \{ (q, p) \in \A :   \frac{1}{q} + \frac{n}{p} \le 1^{-} \} $$

$$ \D = : \{ (q, p ) :  \frac{1}{2q} + \frac{n-1}{4p} \le
(\frac{n-1}{4})^{-} \} $$ 

$$ \D_t = :  \{ (q, p ) : q \ge 2 , \,  \frac{1}{2q} + \frac{n-1}{4p} \le
(\frac{n-1}{4})^{-} \} $$ 

$$\E =: \{ (q, p) : \frac{1}{q} = \frac{1}{q_1} + \frac{1}{q_2}; \,
\frac{1}{p}= \frac{1}{p_1} + \frac{1}{p_2} \text{ with } (q_1, p_1)
\in \A \text{ and } (q_2, p_2)  \in \C \} $$ where $\A$ is as above,
the set of all (wave) admissible pairs. Finally, let 

$$ \G = : \D \cap \E $$

$$ \G_t = :  \D_t \cap \E $$
\endproclaim    

\proclaim{ Remark} Note that  
since  $\A \subset \D$,  $\A \subset \D_t$, and $\A \subset E$
we obviously have that $ \A \subset \G $ and $ \A \subset \G_t $. 
We will refer to the pairs in $\G$ as the set  
of  `good pairs for frequency localized wave products'. 
\endproclaim

\proclaim{(2.10) Definition} Let $\S_{+}^{(-1)}$ be the space of functions on $\IR \times \IR^n$ whose norm 
is given by 
$$ \| \phi \|_{\S_{+}^{(-1)}} \, := \, \sum_{k \in \IZ} \| Q_k (\phi) \|_{{S_k}_{+}^{(-1)}} $$ where 

$$\| \Phi \|_{{\S_k}_{+}^{(-1)}} := \sup_{ (q, p) \in \G} 2^{ k \,
(\frac{1}{q} + \frac{n}{p}-1)}  \| \Phi \|_{L^{q}_t L_x^{p}} \, + \, 
 \sup_{ (q, p) \in \G_t} 2^{ k \,(\frac{1}{q} + \frac{n}{p}-1)} 2^{-k} 
\| \partial_t \Phi \|_{L^{q}_t L_x^{p}}  $$

\endproclaim 

\proclaim{ (2.11) Lemma} We have the following embeddings 

$$\align \S_{+}^{(-1)}  &\hookrightarrow \S^{(-1)} \\
 \S_{+}^{(-1)}  &\hookrightarrow \B_1 \\
\S_{+}^{(-1)}  &\hookrightarrow L^q_t {\dot B}^s_{\tilde p, 2} \qquad  \text{
for all } \, q \ge 2, {\tilde p} \ge 2 \text{ and } s = \frac{1}{q} +
\frac{n}{\tilde p} -1 \endalign $$   

\endproclaim

\demo{Proof} This is an easy consequence of the definition of $\G$, 
Lemma (2.7), the embeddings $\ell^p \subset \ell^q$ for $ p < q$ and
the fact that 
$$\| f \|_{L^q_t {\dot B}^s_{\tilde p, 2}} \lesssim \bigl( \sum_{k \in \IZ} 2^{2
k s} \| Q_k(f) \|^2_{ L^q_t L^{\tilde p}_x}\bigr)^{1/2} $$
We should also note that $(q, {\tilde p} ) \in \G$ for any $q, {\tilde p} \ge 2$ \qquad \qed
\enddemo

\proclaim{(2.12) Proposition (Main Multiplication Estimate)}

$$ P_{-1} : \S^{(-1)} \times \S^{(-1)} \lto  \S_{+}^{(-1)} $$

\endproclaim

\demo{Proof}  
We consider  the first supremum 
term in the $\S_{+}^{(-1)}$-norm; i.e. we need to show : 

$$ \sum_{l \in \IZ} \,  \sup_{ (\tq,\tp) \in \G} 2^{ l (1/\tq + n/\tp -1 )} \| Q_l ( \inv ( f \cdot g) ) \|_{L^{\tq}_t L^{\tp}_x }  \lesssim \|f \|_{S^{(-1)}} \| g \|_{S^{(-1)}} $$

Let $f$ and $g$ be in $\S^{(-1)}$ and let $f_k =Q_k(f)$ and $g_j=Q_j(g)$ be their corresponding Littlewood-Paley projections. We write

$$ \align P_{-1} ( f \cdot g) &=   P_{-1} \bigl(  \sum_{ k, j \in \IZ} f_k \cdot g_j \bigr) \\ 
&=  P_{-1} \bigl(  \sum_{ k, j \in \IZ: k \ge j } f_k \cdot g_j \bigr) + P_{-1} \bigl(\sum_{ k, j \in \IZ : k < j} f_k \cdot g_j \bigr). \endalign  $$ By symmetry of the sums, it is enough to consider only one of them. The proof for the other is identical after exchanging $k$ and $j$

Since $ \text{ supp } \widehat{(f_k \cdot g_{k-m})}\subseteq \{ \xi : |\xi| \le 2^k \} $ we have that $Q_l (f_k \cdot g_{k-m} ) \equiv 0 $ unless $k \ge l$. 
 On the other hand, we have that 
 $ \text{ supp } \widehat{(f_k \cdot g_{k-m})}  \cap \{ \xi : |\xi| << 2^{k-m} \} = \emptyset $ if $ m > 5$
Hence,  $Q_l (f_k \cdot g_{k-m} ) \equiv 0 $ unless $l=k$ and $m > 5$ or $ m \le 5$ and $ l < k $.  

Define ${\tilde Q}_k (f) =  \sum_{ k-5 \le j \le k+5 } Q_j(f)$. By the above argument we conclude that 
it is enough to prove  each of the following two estimates :

$$ \sum_{ l \in \IZ}  \, \sup_{(\tq, \tp) \in \D} \quad
 2^{l (1/\tq + n/\tp -1)} \sum_{ k > l }\| Q_l (\inv ( f_k \cdot {\tilde Q}_k(g)) ) \|_{L^{\tq}_t L^{\tp}_x}  \lesssim \|f \|_{S^{(-1)}} \| g \|_{S^{(-1)}}.\tag2.12)(i$$

$$ \sum_{ l \in \IZ} \,  \sup_{(\tq, \tp) \in \E} \, \,  
2^{-l} 2^{l (1/\tq + n/\tp -1) } \| Q_l\, \bigl( \, \sum_{ m > 5 }  \, f_l \cdot g_{l-m} \bigr) \|_{L^{\tq}_t L^{\tp}_x}  \lesssim  \|f \|_{S^{(-1)}} \| g \|_{S^{(-1)}}.\tag2.12)(ii$$
 
since $\G = \D \cap \E $.

$\bullet$ We consider (2.12)(i). 

\vskip .1in 

For each $(\tq, \tp) \in \D$ let  $s = 1/\tq + n/\tp -1 $

First note that 

$$ \sum_{ l \in \IZ} \, \sup_{(\tq, \tp) \in \D} \, \sum_{ k> l } 
 2^{l s} \| Q_l ( \inv ( f_k \cdot {\tilde Q}_k(g) )) \|_{L^{\tq}_t L^{\tp}_x} \le  \sum_{ k \in \IZ} \,  \sum_{l < k } \, \sup_{(\tq, \tp) \in \D} \,  2^{l (s-1)}  \| Q_l  ( ( f_k \cdot {\tilde Q}_k(g)) \|_{L^{\tq}_t L^{\tp}_x} $$ 

Since $(\tq, \tp) \in \D$ we have that 

$$\frac{1}{2\tq} + \frac{n-1}{4\tp} \le
(\frac{n-1}{4})^{-} .$$ 

By the same argument used in the proof of Lemma (2.7) it is enough to take the supremum over all $(\tq, \tp) \in \D$ such that $\frac{1}{2\tq} + \frac{n-1}{4\tp} =
(\frac{n-1}{4})^{-} .$ We denote this set by $\D^{\#}$. 

Since $\tp < \infty$, let $1 < p < \infty$ such that $$\frac{1}{\tp} < \frac{1}{p} < \frac{1}{\tp} + \frac{1}{50 (n-1)}. $$  Then we have that $$ \frac{n-1}{4 \tp} < \frac{n-1}{4 p} < \frac{n-1}{4 \tp} + \frac{1}{200}. $$ Let $r >1$ be such that $$1 + \frac{1}{\tp} = \frac{1}{r} + \frac{1}{p}.$$ By Young's inequality and H\"older's inequality we then have that 

$$ \align & \sum_{ k \in \IZ} \,  \sum_{l < k } \, \sup_{(\tq, \tp) \in {\D}^{\#}} \,  2^{l (s-1)}  \| Q_l  ( ( f_k \cdot {\tilde Q}_k(g)) \|_{L^{\tq}_t L^{\tp}_x} \\
&\lesssim   \sum_{ k \in \IZ} \,  \sum_{l < k } \, \sup_{(\tq, \tp) \in {\D}^{\#}} \,  2^{l (s-1)}  2^{ n l (1 - 1/r)}  \| f_k \|_{L^{2 \tq}_t L^{2 p}_x} \| {\tilde Q}_k(g) \|_{L^{2\tq}_t L^{2 p}_x}. \endalign $$ 

But by our choice of $p$ we have that if $( \tq, \tp ) \in {\D}^{\#}$ then 
$(2 \tq, 2 p)$ is still in $\A$ the set of admissible pairs. Moreover,

$$ (\frac{n-1}{4})^{-} = \frac{1}{2 \tq} + \frac{ n-1}{ 4\tp } \le \frac{1}{2 \tq} + \frac{ n-1}{ 4p} \le \frac{1}{2 \tq} + \frac{ n-1}{ 4 \tp} + \frac{1}{200} < \frac{ n-1}{4} $$ 

Hence, up to a constant, we can bound the last sum by 

$$ \align &\sum_{ k \in \IZ} \,  \sum_{l < k } \, \sup_{(\tq, \tp) \in {\D}^{\#}} \,  2^{l (s-1)}  2^{ n l (1 - 1/r)}  2^{2 k (1 - \frac{ 1}{ 2 \tq} - \frac{n}{2p})} \| f_k \|_{\S^{(-1)}_k} \| {\tilde Q}_k(g) \|_{\S^{(-1)}_{k}}\\
& \lesssim \sum_{ k \in \IZ} \,  \sum_{j \ge 0 } \, \sup_{(\tq, \tp) \in {\D}^{\#}} \,  2^{-j (\frac{1}{\tq} + \frac{n}{p} - 2)}  \| f_k \|_{\S^{(-1)}_k} \| {\tilde Q}_k(g)  \|_{\S^{(-1)}_{k}}.\endalign $$

Since $n \ge 4$ and $(\tq, \tp) \in  {\D}^{\#}$, we have that $2 \le 2 \tp \le 6$ and hence, 
$$\align \frac{1}{\tq} + \frac{n}{p} \,  > \, \frac{1}{\tq} + \frac{n}{\tp}  &= 2 (\frac{n-1}{4})^{-} - 2  \frac{( n-1)}{ 4 \tp}  + \frac{ 2 n}{2 \tp}  \\
&= 2 (\frac{n-1}{4})^{-}  + \frac{n+1}{ 2 \tp } \\
&>  2 (\frac{n-1}{4})^{-}  + \frac{n+1}{ 6 } \\
&> 2 (\frac{3}{4})^{-} + \frac{5}{6} > 2^{+}.\endalign  $$

Thus, we can sum in $\, j \ge 0\, $ above;  and the desired estimate follows by 
Cauchy-Schwartz in the sum over $k$. 

\vskip .1in 

$\bullet$ We consider (2.12)(ii). 

We proceed as follows,

$$\align &\sum_{ l \in \IZ} \, \sup_{(\tq, \tp) \in {\G}} \, \,  
2^{-l} 2^{l (1/\tq + n/\tp -1) } \| Q_l \, \bigl( \, \sum_{ m > 5 }  \, f_l \cdot g_{l-m} \bigr) \|_{L^{\tq}_t L^{\tp}_x} \\
&\lesssim  \sum_{ l \in \IZ} \, \sum_{ m > 5 } \, \sup_{(\tq, \tp) \in {\E}} \, \,  2^{-l} 2^{l (1/\tq + n/\tp -1) } \| Q_l (  f_l \cdot g_{l-m} ) \|_{L^{\tq}_t L^{\tp}_x} \\
&\lesssim   \sum_{ l \in \IZ} \, \sum_{ m > 5 } \, \sup_{(\tq, \tp) \in {\E}} \, 2^{l (1/\tq + n/\tp -2) }  \| f_l \|_{ L^{q_1}_t L^{p_1}_x} \| g_{l-m} \|_{L^{q_2}_tL^{p_2}} \endalign $$ 
by H\"older's inequality with $\frac{1}{\tq} = \frac{1}{q_1} + \frac{1}{q_2}, \, \,  \frac{1}{\tp} = \frac{1}{p_1} + \frac{1}{p_2}$ and $(q_1, p_1) \in \A$, $(q_2, p_2) \in \C \subset \A$. Recall also that $\| Q_l \|_1 =1$. 

By the Strichartz inequalities we then have that the sum above is bounded by 

$$ \align &\sum_{ l \in \IZ} \, \sum_{ m > 5 } \, \sup_{(\tq, \tp) \in {\E}} \, 2^{l (1/\tq + n/\tp -2) } 2^{ l (2 - 1/\tq - n/\tp )} 2^{- m ( 1 - 1/q_2- n/p_2 )} \| f_l \|_{\S^{(-1)}_l} \| g_{l-m} \|_{ \S^{(-1)}_{l-m}} \\
&\lesssim  \sum_{ m > 5 } \, 2^{- m ( 1/100 ) }  \sum_{ l \in \IZ}  \| f_l \|_{\S^{(-1)}_l} \| g_{l-m} \|_{ \S^{(-1)}_{l-m}} \endalign $$

>from where the desired estimates follows by doing first Cauchy Schwartz in the sum over $l \in Z$ and finally summing over $m > 5$.

\vskip .2in

To obtain the desired estimate for the second supremum in the
definition of the $\S^{(-1)}_{+}$-norm we need to show a companion 
estimates to (2.12) (i) and (2.12) (ii). 

The high-low estimate for a 
time derivative  is treated in a similar manner to (2.12)(ii).
 Indeed, the time derivative first introduces a ``loss'' of $2^l$ and then 
one recoupes $2^l$ (or even $2^{l-m}$) from the estimate for 
$\|\partial_t f_l\|_{L_t^{q_1}L_x^{p_1}}$. We omit the details for
 that part and  we concentrate instead  on the high-high interaction. 
To this end, we will show that 

$$
\sum_l 2^{l(1/\tq+n\tp-3)}\sum_{k\geq l}
\sup_{(\tq,\tp)\in D^{\#}_t}\|Q_l(\partial_tf_k 
\tilde{Q}_k g)\|_{L^\tq_t L^\tp_x}\lesssim \|f\|_{\S_k^{(-1)}} 
\|g\|_{\S_k^{(-1)}},
$$
where $D_t^{\#}=\left\{(\tq,\tp): \tq\geq 2, 
\frac{1}{2\tq}+\frac{(n-1)}{4\tp}=\frac{n-1}{4}-\right\}$.  We note that 
the case, when the time derivative falls on $\tilde{Q}_k g$ is symmetric. 
By applying the same estimates as in (2.12)(i), one obtains that the sum 
 is bounded by 
$$
\sum_l\sum_{l\leq k}\sup_{(\tq,\tp)\in D^{\#}_t}2^{l(1/\tq+n/p-3)} 2^k
2^{2k(1-1/2\tq-n/2p)} \|f\|_{\S_k^{(-1)}} 
\|g\|_{\S_k^{(-1)}},
$$
since by definition 
$\|\partial_t f_k\|_{L^{2\tq}_t L^{2p}_x}\lesssim 2^k 
2^{k(1-1/2\tq-n/2\tp)}\|f\|_{\S_k^{(-1)}}$. 
Therefore, we need to bound 
$$
\sum_l\sum_{l\leq k}
\sup_{(\tq,\tp)\in D^{\#}_t}2^{(l-k)(1/\tq+n/\tp-3)}
\|f\|_{\S_k^{(-1)}} 
\|g\|_{\S_k^{(-1)}}.
$$
which amounts to verifying $1/\tq+n/\tp\geq 3+$, which is somewhat stronger 
than what was  needed in (2.12)(i).  We have 
$$
\frac{1}{\tq}+\frac{n}{p}>\frac{1}{\tq}+\frac{n}{\tp}=
2\left(\frac{1}{2\tq}+\frac{n-1}{4\tp}\right)+\frac{n+1}{2\tp}\geq  
\frac{n-1}{2}+\frac{n+1}{2\tp}-\frac{1}{50}.
$$
The restriction $\tq\geq 2$ in the definition of $D^{\#}_t$ implies 
$\tp\leq (n-1)(n-2)$ and thus 
$$
\frac{1}{\tq}+\frac{n}{p}>\frac{n-1}{2}+\frac{(n+1)(n-2)}{2(n-1)}
-\frac{1}{50}>3 
$$
for $n\geq 4$. We note that one can do a slightly larger domain of 
admissible pairs in the time derivative estimates, but  $D^{\#}_t$ will 
suffice to close the estimates later on.
\enddemo

\vskip .1in

\proclaim{Remark} As a consequence of Theorem (2.12) we have that any element in
$\S^{(-1)}_{+}$ belongs to  both  $\B_1 \subset L^1_t L^{\infty}_x$ and $L^2_t {\dot B}^{n/2 -
1/2}_{2,2}$. This will imply, in particular, that the connection
$1$-form  $a$ - whose existence, uniqueness and regularity is
established in section 3-  belongs to $\B_1 \subset L^1_t L^{\infty}_x$
and $L^2_t {\dot B}^{n/2 - 1/2}_{2, 2}$. This is the crucial fact
needed to obtain the {\it apriori} bounds on the non-linearity
(c.f. Theorem 2.13 below).  

\endproclaim

\proclaim{ Remark} Although the above embeddings will suffice for our purposes in the present paper;  it is interesting to note that in fact, 
$\inv$ maps $ \S^{(-1)} \times \S^{(-1)}$ into a slightly larger class of Besov spaces. Namely into $L^q_t {\dot B}^s_{ \tp, 2}$ for any $ q \ge 2$,  $ s = 1/q + n/\tp -1 $ and $ \tp \ge p$; where  $ p < 2$ is  
such that $\dfrac{1}{q} + \dfrac{n}{p} = \dfrac{n}{2} +1 $. We include a separate proof of this fact in the Appendix. 
\endproclaim
\vskip .1in

\proclaim{(2.13) Theorem} Let $a \in   \S_{+}^{(-1)} $ and $b \in
\S^{(-1)}$ then 

$$ \bigl( \sum_{k \in \IZ}   2^{2k ( n/2 -1 )}  \| Q_k ( a \cdot b )
\|^2_{L^1_t L^2_x} \bigr)^{1/2} \lesssim \| a \|_{\S_{+}^{(-1)}} \| b
\|_{\S^{(-1)}} $$
\endproclaim 

\demo{Proof}
We start as usual by performing a Littlewood-Paley decomposition of $a$ and $b$. We obtain 
$$\align \sum_{k \in \IZ}   2^{2k ( n/2 -1 )}  \| Q_k ( a \cdot b )
\|^2_{L^1_t L^2_x} \, & \lesssim \, \sum_{k \in \IZ}  2^{2k ( n/2 -1 )}  
\| \sum_{m > 5}  Q_k ( Q_k(a) \cdot Q_{k-m}(b) )\|^2_{L^1_t L^2_x} \, \, + \\
& +  \sum_{k \in \IZ}  2^{2k ( n/2 -1 )} \| \sum_{m > 5} Q_k ( Q_{k-m}(a) \cdot Q_{k}(b) ) \|^2_{L^1_t L^2_x} \, \, + \\
& +  \sum_{k \in \IZ} 2^{2k ( n/2 -1 )}  \|\sum_{k < l} Q_k ( Q_l(a) \cdot Q_l(b) )\|^2_{L^1_t L^2_x}.\endalign $$

Now since $a$ and $b$ belong to different spaces we lose the `symmetry' and need to consider all three cases separately. 

$\bullet$  We consider the first of the three sums above. 

$$ \align& \sum_{k \in \IZ} \sum_{m > 5}  2^{ k ( n/2 -1 )}  
\| Q_k ( Q_k(a) \cdot Q_{k-m}(b) )\|_{L^1_t L^2_x} \\
& \lesssim \sum_{m > 5}\sum_{k \in \IZ}   2^{ k ( n/2 -1 )}  
\|   Q_k(a) \|_{L^2_t L^2_x} \|Q_{k-m}(b) \|_{L^2_t L^{\infty}_x} \\
&\lesssim  \sum_{m > 5} 2^{-m/2} \sum_{k \in \IZ}   2^{ k ( n/2 - 1/2 )}  
\| Q_k(a) \|_{L^2_t L^2_x} \|b_{k-m}\|_{\S^{(-1)}_k}, \endalign $$
since the pair $(2, \infty)$ is admissible. Note that the pair $(2, 2) \in \G$, whence $\|   Q_k(a) \|_{L^2_t L^2_x} \le 2^{k (1/2 -n/2)} \| Q_k(a)\|_{{\S_k}_{+}^{(-1)}}$ -since $1/2 + n/2 -1 = n/2-1/2$- . Finally 
 do Cauchy- Schwartz and the desired estimate 
follows after summing over $m > 5$ last.

$\bullet$ We consider next the second sum.

$$\align & \sum_{k \in \IZ}  2^{2  k ( n/2 -1 )}  
\| \sum_{m > 5} Q_k ( Q_k(b) \cdot Q_{k-m}(a) ) \|^2_{L^1_t L^2_x} \\
& \lesssim \sum_{k \in \IZ}  2^{ 2 k ( n/2 -1 )}  
\| Q_k(b) \|^2_{L^{\infty}_t L^2_x} \bigl( \sum_{m > 5}  \| Q_{k-m}(a) \|_{L^1_t L^{\infty}_x}\bigr)^2 \\
&\lesssim \| a \|^2_{{\S}_{+}^{(-1)}} \, \,  \sum_{k \in \IZ} \| Q_k(b) \|^2_{\S^{(-1)}_k}\\
&\lesssim \| a \|^2_{{\S}_{+}^{(-1)}}  \| b \|^2_{{\S}^{(-1)}}\endalign $$

$\bullet$ Finally, we consider the third sum. 

$$\align & \sum_{k \in \IZ} \sum_{k < l }  2^{  k ( n/2 -1 )}  
\| Q_k ( Q_{l}(b) \cdot Q_{l}(a) )\|_{L^1_t L^2_x} \\
& \lesssim \sum_{k \in \IZ} \sum_{k < l }  2^{  k ( n/2 -1 )}  
\| Q_{l}(b) \|_{L^{\infty}_t L^2_x} \| Q_{l}(a) \|_{L^1_t L^{\infty}_x}\\
& \lesssim \sum_{l \in \IZ} \sum_{j \ge 0 }  2^{ -j ( n/2 -1 )}  
 \| Q_{l}(b) \|_{\S^{(-1)}_l} \| Q_{l}(a) \|_{L^1_t L^{\infty}_x} \\
& \lesssim \sum_{l \in \IZ}  \| Q_{l}(b) \|_{\S^{(-1)}_l} \| Q_{l}(a) \|_{L^1_t L^{\infty}_x}, \endalign $$
>from where by Cauchy-Schwartz we obtain the desired estimate invoking once again the fact that $a \in \B_2$ \qquad \qed
\enddemo
\vskip .5in

\head{3. Existence, uniqueness and regularity of the connection 1-form 
in $\S^{(-1)}$ }\endhead

\vskip .1in 

\proclaim{ Proposition 3.1} Let $b \in \S^{(-1)}$ have sufficiently small norm; then the map $$ \Phi (w) = \inv [ w, w] + \inv [b, b]$$ has a unique fixed point $a = \Phi (a) 
\in \S^{(-1)}$. Moreover, the fixed point $a$ belongs to $ L^1_t L^{\infty}_x \cap {\dot B}^s_{\tp,2} $ for any $ \tp \ge 2$ and $  n/\tp - 1 \le s \le n/\tp - 1/2$ .\endproclaim 

\demo{Proof} Let $b \in \S^{(-1)}$ such that $\| b \|_{ \S^{(-1)}} = \epsilon \le \dfrac{1}{100 c_1}$ where $ c_1 >0$ is the constant from (2.12) such that 
$\| \inv [ \, \alpha, \beta\, ] \|_{\S^{(-1)}_{+}} \le c_1 \| \alpha \||_{\S^{(-1)}} \| \beta\|_{\S^{(-1)}} $.

Let $0< r < \epsilon $ and let $B=B_r(0) $ be the ball in $\S^{(-1)}$ centered at $0$ and radius $r$. Then 

\roster
\item $ \Phi: B_r (0) \to   B_r (0) $

\item $ \| \Phi(w_1) - \Phi(w_2) \|_{\S^{(-1)}} <   \| w_1  - w_2 \|_{\S^{(-1)}} $

\endroster 

To check (1) let $w \in B_r(0)$ then $$\| \Phi(w)\|_{\S^{(-1)}} \le c_1 \| w \|^2_{\S^{(-1)}} + c_1 \| b \|^2_{\S^{(-1)}} < r $$
To check (2) let $w_1, w_2 \in  B_r(0)$ then $$\align \| \Phi(w_1) -\Phi(w_2) \|_{\S^{(-1)}} &\le c_1 \| w_1 - w_2 \|_{\S^{(-1)}} \max_{i =1, 2} \{ \| w_i\|_{\S^{(-1)}} \}\\
& \le \frac{1}{100} \| w_1 - w_2 \|_{\S^{(-1)}}  < \| w_1  - w_2 \|_{\S^{(-1)}} \endalign $$

Thus $\Phi$ is a contraction and  hence there exists a unique fixed point $a = \Phi(a) \in \S^{(-1)}$ such that 
$$ a = \inv [a, a] + \inv [b, b] $$ 

By Lemma 2.12,  $\, \inv [\cdot , \cdot ] \in L^1_t L^{\infty}_x \cap {\dot B}_{\tp, 2}^s$ for any $ \tp \ge 2$ and $ s= 1/\tq + n/\tp -1$ with $\tq \ge 2$ . Hence so does $a$ \qed. 
\enddemo

\head{4. The Modified Wave Map System }\endhead

\vskip .1in

In this section we prove that the Cauchy problem for the MWM system derived in Section 2 has a unique global solution in $L^{\infty}(\IR; {\dot H}_x^{n/2} ) $ provided the initial data has sufficiently small ${\dot H}_x^{n/2} \times {\dot H}_x^{n/2-1}$ norm.  

Let us denote by $B(a, b)$ the quadratic form equal to any finite linear combination of functions $a \in \S^{(-1)}_{+} $ and  $b \in \S^{(-1)} $ of the form $\sum_{\kappa, \ell} c_{\kappa \ell}  a_{\kappa} b_{\ell}$ where $a_{\kappa} \in \S^{(-1)}_{+} $,  $b_{\ell} \in \S^{(-1)} $ and $c_{\kappa \ell} \in \IC$. 

According to our reductions in the previous section, we consider the system 
of coupled wave equations in $R^{n+1}$, $ n \ge 4$. 

$$\align \square v &= B(a, b)\\
v(x, 0) &= f(x) \tag4.1\\
v_t(x, 0) &= g(x) \\
\endalign  $$

\proclaim{(4.2) Lemma } Let $ a \in \S^{(-1)}_{+}$ and $b \in \S^{(-1)} $. Then the solution to the MWM system (4.1) with initial data $(f, g) \in {\dot H}^{n/2} \times {\dot H}^{n/2-1}$ satisfies   
$$ \| v \|_{\S} \lesssim \| f \|_{ {\dot H}^{n/2}} + \| g \|_{{\dot H}^{n/2-1}} + \| a \|_{\S^{(-1)}_{+}} \| b \|_{\S^{(-1)}} $$ 
\endproclaim

\demo{Proof} Let us denote by $v_k = Q_k(v)$. By the Strichartz's estimates, we have that 

$$ \| v_k \|_{\S_k} \lesssim  \| f_k \|_{ {\dot H}^{n/2}} + \| g_k \|_{{\dot H}^{n/2-1}}  + 2^{k ( n/2 -1) } \| B( a, b) \|_{L^1_t L^2_x} $$ from where by Theorem (
2.13) we have that 

$$ \align  \| v \|_{\S} &  = \bigl( \sum_{k \in \IZ} \| v_k \|_{\S_k}^2 \bigr)^{1/2} \\
&\lesssim \| f \|_{ {\dot H}^{n/2}} + \| g \|_{{\dot H}^{n/2-1}} + \bigl( \sum_{k \in \IZ}  2^{ 2 k ( n/2 -1) } \| B( a, b) \|^2_{L^1_t L^2_x} \bigr)^{1/2} \\
&\lesssim \| f \|_{ {\dot H}^{n/2}} + \| g \|_{{\dot H}^{n/2-1}} + \| a \|_{\S^{(-1)}_{+}} \| b \|_{\S^{(-1)}}  \endalign $$ as desired \qquad \qquad \qquad \qquad \qed . 

\enddemo 

\proclaim{ (4.3) Theorem (Existence)} 

There exists $\varepsilon > 0$ such that whenever the initial data $\| (f, g ) \|_{{\dot H}^{n/2} \times {\dot H}^{n/2-1}} < \varepsilon $, 
the system (4.1) has a unique global solution $v \in \S $.
 
In particular, the solution $v$ belongs both to 
$$\bullet \,  L^{\infty}( \IR; {\dot H}_x^{n/2} ) \, \cap \, L^2( \IR; {\dot B}^{1}_{2n, 2}) \qquad {\text {and}} $$  
$$ \bullet \,  W^{1, {\infty}}( \IR; {\dot H}_x^{n/2-1} ) \,\cap \, W^{1,2}( \IR; {\dot B}^{0}_{2n, 2}).$$  
Moreover, there is  {\it stability; i.e.} 
$$  {{\text{\rm ess}}}\sup_{t} \| v_1 - v_2 \|_{ {\dot H}^{n/2}} \lesssim \| (f_1, g_1) - (f_2, g_2) \|_{{\dot H}^{n/2} \times {\dot H}^{n/2-1}} $$   provided the r.h.s. is small enough.
\endproclaim

\demo{Proof}
The proof proceeds by  Picard's iteration relying on the {\it a priori} estimates as well as the necessary smalless of the data.

Suppose  $\| ( f, g)\|_{{\dot H}^{n/2} \times {\dot H}^{n/2-1}}= \delta $ and let $v_0$ be the solution to 
$$\square v_0 = 0; \qquad v_0(0, \cdot) = f \quad \partial_t v_0(0, \cdot) = g .$$
By the Strichartz's estimates
$$ \| v_0 \|_{\S} \le c_1 \| ( f, g)\|_{{\dot H}^{n/2} \times {\dot H}^{n/2-1}} = c_1 \delta .$$

Now, $v_0 = ( \varphi_0, \psi_0)$ produces $b_0 = d \, \varphi_0 + \, div_{(sp.t)} \, \psi_0 $ with $\| b_0\|_{\S^{(-1)}} \le c_2 \| v_0 \|_{\S} \le c_3 \delta $. 

Next, the multiplication estimates allow one to perform a 
 fixed point argument to produce $a_0$ from $b_0$ by solving 
$$ a_0 = \inv [ a_0, a_0] + \inv [b_0, b_0]. $$ Moreover, 
$\|a_0\|_{\S^{(-1)}_{+}} \le c_4 \|b_0\|_{\S^{(-1)}}^2 \le c_5 {\delta}^2$ 

Let $v_1$ be the solution of 
$$\square v_1 = B(a_0, b_0) \qquad v_1(0, \cdot) = f \quad \partial_t v_1 (0, \cdot) = g .$$ By the {\it a priori} estimate,  
$$ \| v_1\|_{\S} \le c_0 \bigl( \delta + \|a_0\|_{\S^{(-1)}_{+}}  \|b_0\|_{\S^{(-1)}} \bigr) \le  2 c_0  \delta $$ provided $\delta$ is small enough.

We proceed next by induction to show that for any $j \ge 0$,  $\Vert b_j \Vert_{\S} \le 2 c_2 c_o \delta$, \,  $\Vert a_j \Vert_{\S} \le c_5\delta^2$ and thus 
$\|v_{j +1} \|_{\S} \leq 2 c_0 \delta$ provided $\delta > 0$ is small enough 
(indep. of $j$),  where 
$v_{j+1}$ is the solution to 
$$\square v_{j+1}  = B(a_j, b_j) \qquad v_{j+1} (0, \cdot) = f \quad \partial_t v_{j+1}(0, \cdot) = g .$$ Note that once again by the {\it a priori} estimates 
$$ \| v_{j+1}\|_{\S} \le c_0  \bigl( \| ( f, g)\|_{{\dot H}^{n/2} \times {\dot H}^{n/2-1}}  + \|a_j\|_{\S^{(-1)}_{+}}  \|b_j\|_{\S^{(-1)}} \bigr). $$

Lastly, for the differences,
$$\align  \square (v_{j+2} - v_{j+1}) &= B(a_{j+1}, b_{j+1}) - B(a_j, b_j)= B( a_{j+1}-a_j, b_{j+1})+ B(a_j, b_{j+1}- b_j) \\
 v_{j+1}(0, \cdot\,) &= 0 \\
\partial_t v_{j+1} (0, \cdot\,) &= 0 .\endalign$$ On the other hand note that since 
$$ a_{j+1} - a_j  = \inv [ a_{j+1} - a_j , a_{j+1}] + \inv [a_{j},  
 a_{j+1} - a_j ] + \inv [b_{j+1} - b_j  , b_{j+1}] + \inv [ b_{j} ,   
b_{j+1} - b_j ] $$
 
$$\| b_{j+1} - b_{j} \|_{\S^{(-1)}} \le c_2 \| v_{j+1} - v_{j} \|_{\S}$$ and 

$$ \|a_{j+1} \|_{\S^{(-1)}}\, , \|a_{j} \|_{\S^{(-1)}}\, , \|b_{j+1} \|_{\S^{(-1)}}\, , \|b_{j} \|_{\S^{(-1)}}\, \le c \delta\, ,  $$
 we have that 
$$\align  \| a_{j+1} - a_j\|_{\S^{(-1)}_{+}} &\le  c \|a_{j+1} \|_{\S^{(-1)}} \| a_{j+1} - a_j \|_{\S^{(-1)}} + \|a_{j} \|_{\S^{(-1)}} \| a_{j+1} - a_j \|_{\S^{(-1)}} + \\ 
&+ \|b_{j+1} \|_{\S^{(-1)}} \| b_{j+1} - b_j \|_{\S^{(-1)}} + \|b_{j} \|_{\S^{(-1)}} \| b_{j+1} - b_j \|_{\S^{(-1)}}\\
&\le c \delta \| a_{j+1} - a_j \|_{\S^{(-1)}} + c' \delta  \| b_{j+1} - b_j \|_{\S^{(-1)}}.\endalign $$ Hence, 
$$ \| a_{j+1} - a_j\|_{\S^{(-1)}_{+}} \le c \delta \| b_{j+1} - b_j \|_{\S^{(-1)}} \le c \delta  \| v_{j+1} - v_{j} \|_{\S}.$$
All in all we then have that 
$$\align \| v_{j+2} - v_{j+1} \|_{\S} &\le c \bigl( \| a_{j+1} - a_j\|_{\S^{(-1)}_{+}}  \|b_{j+1} \|_{\S^{(-1)}} +  \|a_{j+1} \|_{\S^{(-1)}} \| b_{j+1} - b_{j} \|_{\S^{(-1)}} \bigr) \\
&\le c \delta^2  \| v_{j+1} - v_{j} \|_{\S} \endalign $$
Finally, by choosing $\delta$ small enough we have that 
$$ \| v_{j+2} - v_{j+1} \|_{\S} \le \frac{1}{2} \| v_{j+1} - v_{j} \|_{\S}. $$

Hence $v_j$ is Cauchy in $\S$, thus establishing existence and uniqueness. 
For the stability result one proceeds in the same fashion as in the proof of being Cauchy; thus concluding the proof of the theorem.  \qed.  

\enddemo

The theorem above gives uniqueness solely in $\S$ which is not enough to claim the solution to the MWM system came from a wave map. Thus we proceed next
to prove a stronger uniqueness result which will indeed suffice in section 5 to return to the wave map.

\proclaim{ (4.4) Theorem (Uniqueness)} Suppose $(v_1\, , \, a_1)$ and $(v_2\, , \, a_2)$ 
are two solutions to 
$$ \align &\square v + B( a, d v) \,= \,0\\
&\Delta a + div B(a, a) + div B(dv, dv )\, =\, 0 .\endalign $$ such that $dv_j = b_j$, for $j=1,2$ are small in $L^{\infty}_t L^n_x$. Suppose that $dv_j = b_j \in L^2_t L^{2n}_x$ for $j=1,2$. Assume in addition that  $a_1 = a_1(v_1) \in L^1_t L^{\infty}_x$. Then $v_1 \, = \, v_2$.   
\endproclaim 

\proclaim{Remarks} The smallness of $d v_j $ in $L^{\infty}_t L^n_x$ is the necessary condition to solve the `gauged' equation. Also note that it is \underbar{not} necessary for $a_2 = a_2 (v_2) \in L^1_t L^{\infty}_x$.
\endproclaim

\demo{Proof} The proof follows the scheme devised by Shatah-Struwe to establish
uniqueness \cite{\SHSb}, \cite{\SHS}. Let us denote  $$ \align \delta w &=  v_1 -  v_2 \\
\delta a &= a_1 - a_2 \\
\delta b &= d v_1 - d v_2 \, , \endalign $$ and so on. Then, 
$$ \align \square \delta w &= B( a_1 , \delta b) + B( \delta a, b_2 )\\
\int \langle \, \frac{\partial}{\partial t} \delta w\, &, \, \square \delta w \, \rangle \, dx \, = \, \frac{1}{2}\, \frac{\partial}{\partial t} E^2 \, , \quad \text{ where }\\
E^2 \,= \, \int | d (\delta w) |^2 &+ |\frac{\partial (\delta w)}{\partial t}|^2 \, dx \, = \, \int |\delta b|^2 \, dx.  \endalign $$ Then, 
$$\frac{1}{2}\, \frac{\partial}{\partial t} E^2 \le \Vert a_1 \Vert_{t, L^{\infty}_x} E^2(t) \, + \, E(t) \Vert \delta a\Vert_{t, L^{\frac{2n}{n-1}}_x} \Vert b \Vert_{t, L^{2n}_x}. $$ Integrating over $t$ we then obtain that 
$$\align E^2(t) &\le \max_{ \tau \le t} E^2(\tau) \, \int_0^{\tau} \Vert a_1 \Vert_{t, L^{\infty}_x} \, dt \\
&+  \max_{ \tau \le t} E(\tau) \, \Vert \delta a\Vert_{L^2_{(0,t)} L^{\frac{2n}{n-1}}_x} \Vert b \Vert_{L^2_{(0,t)} L^{2n}_x} \, ,  \tag4.5 \endalign$$ where $L^2_{(0,t)}$ means the $L^2$ norm on the time interval $(0, t)$. Now, 
$$\Delta \delta a + div B(\delta a, a_1 + a_2) + div B( \delta b, b_1, b_2 ) = 0$$ Hence, $$\align & \Vert \delta a\Vert_{(t, W_x^{1, 2n/n+1})} \le \\
& \le\Vert \delta a \Vert_{(t, L^{2n/n-1}_x)} ( \Vert a_1 \Vert_{(t, L^n_x)} + \Vert a_2 \Vert_{(t, L^n_x)} ) \, + \, \Vert \delta b \Vert_{(t, L^{2}_x)} ( \Vert b_1 \Vert_{(t, L^{2n}_x)} + \Vert b_2 \Vert_{(t, L^{2n}_x)} ). \tag4.6 \endalign$$ On the other hand, $\Vert a_j\Vert_{W^{1, n/2}} \le \Vert a_j \Vert^2_{L^n} + \Vert b_j \Vert^2_{L^n}  $ and $\Vert b_j \Vert^2_{L^n}$ is small for each $t$. Moreover, by Sobolev embedding 
$\Vert a_j \Vert_{L^n} \le c(n) \Vert a_j \Vert_{{\dot W}^{1, n/2}} $; hence (for example by a fixed point argument in $L^n$ similar to Lemma 3.1 ) we have that 
$\Vert a_j \Vert_{L^n}$ is also small for each fixed $t$. 

All in all, from (4.6) we have that, 

$$ \Vert \delta a \Vert_{t, L^{\frac{2n}{n-1}}_x} \le {\tilde c}(n) \Vert \delta a \Vert_{t, L^{\frac{2n}{n+1}}_x} \le {\tilde c}(n) E(t) ( \Vert b_1 \Vert_{(t, L^{2n}_x)} + \Vert b_2 \Vert_{(t, L^{2n}_x)} ). $$ 

Integrate $\tau \le t$ to get,

$$ \Vert \delta a \Vert_{L^2_{(0,t)}, L^{\frac{2n}{n-1}}_x} \le {\tilde c}(n) \max_{\tau \le t} E(\tau) ( \Vert b_1 \Vert_{L^2_{(0,t)}, L^{2n}_x} + \Vert b_2 \Vert_{L^2_{(0,t)}, L^{2n}_x} ). $$

Sticking this estimate back in (4.5) we obtain 

$$ E^2(t) \le \max_{\tau \le t} E^2 (\tau) \bigl( \Vert a_1\Vert_{L^1_{(0, t)} L^{\infty}_x} \, + \, {\tilde c}(n) ( \Vert b_1 \Vert_{L^2_{(0,t)}, L^{2n}_x} + \Vert b_2 \Vert_{L^2_{(0,t)}, L^{2n}_x} )^2 \bigr).$$ Since $E(0)=0$, we must then have that $E(t) =0$. \qed. 
\enddemo

By differentiating the MWM system (4.1) and observing that the
resulting nonlinearity has the same bilinear structure 
-for which the main multiplication estimates
hold- the following regularity result follows. 
 
\proclaim{ (4.7) Theorem ( Higher Regularity )} Suppose the initial
data $(f, g)$ to (4.1) is in $H^{n/2 +1} \times H^{n/2}$ and has
sufficiently small ${\dot H}^{n/2} \times {\dot H}^{n/2-1}$ norm. Then
the solution $v$ to the Cauchy problem (4.1) with initial data $(f,
g)$ can be continued in  $H^{n/2 +1} \times H^{n/2}$ globally in
time. Furthermore, we have the global bounds $$ \Vert v 
\Vert_{L_t^{\infty}(\IR; {\dot H}_x^{n/2 +1} )}
\lesssim \Vert (f, g ) \Vert_{ {\dot H}_x^{n/2 +1} \times {\dot
H}_x^{n/2}}. $$ 
\endproclaim
\demo{Proof}
Assume for simplicity that the data is infinitely smooth. The constants in 
our estimates will depend only on the relevant smoothness assumptions in the 
theorem.

Differentiate  (4.1) to get
$$\align \square w &= B(\partial a, b)+B( a,\partial b)\\
w(x, 0) &= \partial f(x)\\
w_t(x, 0) &= \partial g(x), \\
\endalign  $$
where $w=\partial v$ ($\partial$ may signify any of 
$\partial_j$). Recall also 
that $a$ is a (unique) fixed point for $\Phi$ and therefore its derivative 
will satisfy 
$$
\partial a=P_{-1} [\partial a,a]+  P_{-1}[a, \partial a] + 
P_{-1} [\partial b,b]+ P_{-1}[b, \partial b].
$$
Estimating both sides in $\|\cdot\|_{\S_{+}^{-1}}$, together with the 
main multiplication estimate and 
$\|\cdot\|_{\S_{+}^{-1}}\leq \|\cdot\|_{\S^{-1}}$, yields
$$
\|\partial a\|_{\S_{+}^{-1}}\lesssim 
\|\partial a\|_{\S_{+}^{-1}}\|a\|_{\S^{-1}}+
\|\partial b\|_{\S^{-1}}\|b\|_{\S^{-1}}.
$$
Recall from the Pickard iteration method, that since 
$\|(f,g)\|_{\dot{H}^{n/2}\times \dot{H}^{n/2-1}}$ is small, we have 
$\|a\|_{\S_{+}^{-1}}$ and $\|b\|_{\S^{-1}}$  small as well. By the usual hiding argument, 
one deduces 
$$
\|\partial a\|_{\S_{+}^{-1}}\lesssim 
\frac{\|\partial b\|_{\S^{-1}}\|b\|_{\S^{-1}}}{1-\|a\|_{\S^{-1}}}\leq
\| b\|_{\S}\|b\|_{\S^{-1}},
$$
thus placing the nonlinearity $B(\partial a,b)$ in the form 
$\S^{-1}_{+}\cdot \S^{-1}$ as in Lemma 4.2. Since 
$b\sim \partial v$, the same holds for the other 
nonlinearity  associated with the derivated equation, namely 
$B(a,\partial b)$. An application of the Strichartz 
estimates and  Lemma 4.2 yields 
$$\align
\|w\|_{\S} & \lesssim \|\partial f\|_{\dot{H}^{n/2}}+ 
\|\partial g\|_{\dot{H}^{n/2-1}}+\|b\|_{\S^{-1}}^2 \| b\|_{\S}\lesssim  \\
& \lesssim  \|(f,g)\|_{H^{n/2+1}\times H^{n/2}}+
o(\|(f,g)\|^2_{\dot{H}^{n/2}\times \dot{H}^{n/2}})\|w\|_{\S}.
\endalign
$$
The result follows, since $\|v\|_{L^{\infty}_t H_x^{n/2+1}}\lesssim 
\|w\|_{\S}. $
\enddemo
\vskip .5in
\head{5. The Return to the Map   }\endhead

The well-posedness results on the modified wave map apply to a larger class of formal solutions $(a, b)$ to the equation than those which come from wave maps. Our method of using the results on the modified wave map equation to show existence of wave maps is similar to the idea we used for non-linear Schr\"odinger \cite{\NSU} and not very different from the technique used by Shatah-Struwe \cite{\SHS}. 
The translation depends on the compactness of $M$ (or certain bounds on the isometric Nash embedding of a non-compact $M$ in an Euclidean space). The proofs are very simple for the Lie group case because of the natural parallel structure; and the compact symmetric space case (e.g. $\IS^m$) is a special case due to the totally geodesic embedding $G/K \subseteq G$. Since we have estimates only for this case, we restrict to this case; although the theorems below are 
true in general. 

\proclaim{ (5.1) Theorem} Let $n \ge 3$. If $( s, \, s^{-1} \, \dfrac{\partial s}{\partial t} ) \in {\dot H}^{n/2} \times {\dot H}^{n/2-1}$ are sufficiently small initial data for a wave map into a compact Lie group $G$, then there exists a gauge transformation $g \in {\dot H}^{n/2}$ and a formal derivative $\dfrac{\partial g}{\partial t} \in {\dot H}^{n/2 -1}$, such that the initial data $$ b = \frac{1}{2} g ( s^{-1} d s) g^{-1} $$ are small in ${\dot H}^{n/2-1}$. Furthermore, if $$a= - dg\, g^{-1} + \frac{1}{2} g ( s^{-1} ds) g^{-1}, $$ then \,$a$\, satisfies $\sum_{j=1}^n \dfrac{\partial a_j}{\partial x^j} =0 $ and is small in ${\dot W}^{n/2, 2n/(n+2)} \, \subset \, {\dot H}^{n/2 -1 } $

\endproclaim

\demo{Proof} Note that the pull-back connection in the frame of left pull-back to the Lie algrebra is $d + \dfrac{1}{2} s^{-1}  ds $. The curvature is $\dfrac{1}{4} [\, s^{-1} ds\, , \, s^{-1} ds\, ]$, which will be small in $L^{n/2}$ since ${\dot H}^{n/2 -1} \subset L^n .$ We can then apply Theorem (1.1) (or actually the first step in a time-slice of the proof) to get a good gauge. Since 
$$ \sum_{j=1}^n \frac{\partial}{\partial x^j} \bigl[ (d g \, g^{-1}) - g ( \frac{ s^{-1} d s }{2} ) g{-1} \bigr] \, = \, 0 , $$ a standard regularity theorem will give $g$ to be as smooth as $s$. Here we use heavily the fact that $g$ is bounded. Then $b = g ( s^{-1} ds ) g^{-1} $ has components which are small in ${\dot H}^{n/2-1}$. 
Since $a = - dg \, g^{-1} + g (\frac{s^{-1} d s}{2}) g^{-1} $ has $\sum_{j=1}^n \frac{\partial a_j}{\partial x^j} \, = \, 0$ on the time-slice $t =0$; $a \in {\dot W}^{1, n/2}$ is small. A standard regularity theorem applied to the equation 
$$ \Delta a_j + \sum_{k=1}^n \frac{\partial}{\partial x^k} [ a_k, a_j ] +   \sum_{k=1}^n \frac{\partial}{\partial x^k} [ b_k, b_j ] \, =\,0 $$ gives $a_j \in {\dot W}^{n/2, 2n/(n+2)}$ small and bounded by $\Vert b \Vert_{{\dot H}^{n/2-1}}$. 

The time derivative $\dfrac{\partial g}{\partial t}$ is chosen so that if $$ \align a_0 & = - \frac{\partial g}{\partial t} \, g^{-1} + \frac{1}{2} g^{-1} s^{-1} \frac{\partial s}{\partial t} g\, , \\
0 &= \Delta a_0 + \sum_{j=1}^n \, \frac{\partial}{\partial x^j} \bigl( [ a_j , a_0] + [b_j, b_0] \bigr)\, , \endalign $$ and $a_0$ will also be small in ${\dot W}^{n/2, 2n/(n+2)}$. This implies $g^{-1} \dfrac{\partial g}{\partial t}$ is small in ${\dot H}^{n/2 -1}$ as claimed. 

The estimates follow from standard composition and multiplication theorems, and elliptic regularity. The needed multiplication theorems are less straightforward for the fractional derivatives needed in odd dimensions, but are extended to the fractional derivatives by interpolation between integral derivatives.  \qed

\enddemo

\proclaim{ (5.2) Theorem} Let $( s, s^{-1} \frac{\partial s}{\partial t} ) \in {\dot H}^{n/2} \times {\dot H}^{n/2 -1}$ be initial data for a wave map into a compact group. If $s$ is sufficiently close to the identity in ${\dot H}^{n/2}$, then there exist approximations $( s_{\alpha}, \nu_{\alpha}) $ 
in $C^{\infty}_{id} \times C^{\infty}$ such that 
$$( s_{\alpha}, \nu_{\alpha}) \to ( s, s^{-1} \frac{\partial s}{\partial t} ) \quad \text{ in } \quad {\dot H}^{n/2} \times {\dot H}^{n/2 -1}.$$

\endproclaim 
By $C^{\infty}_{id}$ we have denoted the space of $C^{\infty}$ maps which are the identity  at infinity. 
\demo{Proof} Let $s \in G \subset \IR^{\ell} \times \IR^{\ell}$, and $ \nu= s^{-1} \frac{\partial s}{\partial t} \in \gG $, which is a vector space. The standard approximation method is to convolve  
$$ \nu_{\alpha} (x) = \int \, \nu( x + 2^{-\alpha} y ) \varphi (y) \, dy = \int \, \nu( x + y' ) \varphi_{\alpha} (y') \, dy' = J_{\alpha} (\nu) (x) $$ where $ \varphi_{\alpha}( y') = 2^{n \alpha} \varphi ( 2^{\alpha} y')$ and $\varphi$ is a smooth bump function with compact support such that $ \int \varphi =1 $. Since $\nu = s^{-1} \frac{\partial s}{\partial t}$ is in the Lie algebra, this makes sense. The approximation for $s$ is more subtle. Let $ \IP : \U(G) \to G$ be the projection operator of a neighborhood of $G$ onto the nearest point in $G$. We define 
$$ s_{\alpha} = \IP( J_{\alpha} (s) ). $$  
This is well define in the case that $$\int_{|x-y| \le r} \, | d s| \,
\le \, r^{-n+1} \, \epsilon$$ for all small $r>0$ and $\epsilon>0$
sufficiently small (depending on the diameter of the neighboorhood
$\U(G)$ ). If $s \in {\dot H}^{n/2}$ is sufficiently small, this will
be true. Then, the result that $s_{\alpha} \to s$ in ${\dot H}^{n/2}$
follows by applying the regularity or density result of F. Bethuel 
of smooth maps between certain manifolds in Sobolev spaces.   
( \cite{\BET} \cite{\BETb} and references therein) . \qed

\enddemo

\proclaim{ (5.3) Theorem} Let $s : [0,T] \times \IR^n \to G$ be a wave map in a time interval $[0,T]$ such that $d s \in L^{\infty}_t L^{n}_x \cap 
 L^2_t  L^{2n}_x$. Assume the initial data is in ${\dot H}^{n/2} \times {\dot H}^{n/2-1}$ and has sufficiently small norm. Then $s$ is a gauge transformation of a modified wave map, and $s \in {\dot H}^{n/2} \times {\dot H}^{n/2 -1 }$ remains small in $[0,T]$. Moreover, if the initial data is in ${\dot H}^{n/2+1} \times {\dot H}^{n/2}$, then $ s \in L^{\infty}_t {\dot H}^{n/2 +1} $ and $ \frac{\partial s}{\partial t } \in L^{\infty}_t {\dot H}^{n/2} $ for the time the solution exists. 
\endproclaim

\demo{Proof} Since solutions of the wave map are local, we can assume without loss of generality that its $L^{\infty}_t L^n_x$ norm is small for the time interval of existence ( {\it a posteriori} this will be true anyway ). Make a gauge transformation to a modified wave map. The gauged modified wave map lies in the regime of our uniqueness theorem (4.4). Therefore, it coincides with the solution we have found (the constructed solution satisfies $ a \in L^1_t L^{\infty}_x $). Hence it is a gauge transformation of a solution in $\S$. 

The regularity theorem (4.7) implies the second statement. \qed.
\enddemo
We define next $\tilde S$ as the natural mixed Lebesgue normed space $S$ lies in. More precisely, 
 
\proclaim{ (5.4) Definition} Let $\tilde S$ be the space of functions on $\IR \times \IR^n$ whose norm is given by 
$$ \| \phi \|_{\tilde S} := \bigcup_{(q, p, s) : \frac{1}{q} + \frac{n}{p} = s} \|\phi \|_{ L^{q}_t {\dot W}^{s, p}}. $$   
\endproclaim 

\proclaim{ (5.5) Corollary} Suppose $s : [0,T] \times \IR^n \to G$ is a wave map with $d s \in L^{\infty}_t L^n_x \cap L^2_t L^{2 n}_x $. Suppose, in addition, the data at any point of time is small in ${\dot H}^{n/2} \times {\dot H}^{n/2 -1}$. Then $\,s\, $ exists for all time and $d s \in {\tilde S} $. 
\endproclaim 

\demo{Proof} The gauge transformation of this map coincides with the MWM we have found. Moreover, if $a\, b, \in S$, then solution $g_{\pm}$ of 
$$ d g + a q \pm b g = 0 $$ exist -since the curvature of $  d + a \pm b$ is zero, we can apply Theorem 1.1 - and a standard regularity argument shows that $d g_{\pm} \in \tilde S$. Then $s = g_{+} g^{-1}_{-}$ has the same property. 

\enddemo 

\proclaim{ (5.6) Theorem } If $(s, \nu) \in H^{n/2+1} \times H^{n/2}$ are initial data for a wave map and $(s, \nu) \in {\dot H}^{n/2} \times {\dot H}^{n/2-1}$ has small enough norm, then there exists a unique global solution with $d^2 s \in {\tilde S}$.
\endproclaim

\demo{Proof} Local existence theorems for data in $H^{n/2+1} \times H^{n/2}$ are available (\cite{\KLM} \cite{\KLS}). By theorem 5.3, the norm of $( s, s^{-1} \frac{\partial s}{\partial t}) \in L^{\infty}_t H^{n/2+1}_x \times  L^{\infty}_t H^{n/2}_x $ remains bounded. Hence the local existence theorems can be used to extend the solution intervals to obtain a unique global solution. 
\enddemo

\proclaim{ (5.6) Theorem } Let $(s, \nu) \in {\dot H}^{n/2} \times {\dot H}^{n/2 -1}$ be small data for a wave map into a compact group or symmetric space. Then there exists a global solution, which is a gauge transformation of a solution to the modified wave equation in $\S$ and hence $ds \in \tilde S$. 
\endproclaim

\demo{Proof} Approximate $(s, \nu)$ by smooth data $(s_{\alpha}, \nu_{\alpha})$. Then there exist global solutions to the wave map problem with initial data $(s_{\alpha}, \nu_{\alpha})$. These are gauge transformations of solutions of the modified wave map problem. Choose a weak limit. This limit must have a gauge transformation coinciding with one of our constructed solutions to the modified wave map problem. But since this solution is a weak limit of solutions satisfying $ d a + [ a, a] + [b, b] = 0$, this modified wave map has properties of the complete wave map and can be gauged back (using theorem (1.4) since the curvature of $d + a \pm b $ is zero). \qed.  

\enddemo

\vskip .3in
\head{Appendix}\endhead

We provide in this section alternate proofs of various multiplication Lemmas which are included in the Main Multiplication Lemma. We start with an auxiliary Lemma.

\proclaim{ (A.1) Lemma} Let $n \ge 4$ and let $f$ be a function on $ \S^{(-1)}$. For any  $q \ge 2$, and
$p^{\prime}$ defined by $ \frac{1}{q} + \frac{n}{p'} =1 $ we have that
$$ \bigl( \sum_{k \in \IZ} \|Q_k(f)\|_{L^q_t
L^{p'}_x}^2 \bigr)^{1/2} \, \lesssim \| f \|_{\S^{(-1)}} \tag A.1)(i$$ 
In addition, we also have 
$$ \|  f   \|_{L_t^\infty \dot{W}^{n/2-1, 2}} \lesssim  \bigl( \sum_{k \in \IZ} 2^{2k (n/2 -1)} \|Q_k(f)\|_{L^\infty_t
L^2_x}^2 \bigr)^{1/2} \, \lesssim \| f \|_{\S^{(-1)}} \tag A.1)(ii$$
The same conclusions hold for $2^{-k} \| \partial_t Q_k(f) \| $ replacing $\|Q_k(f)\| \, \,$ (in the corresponding norms) . 
 
\endproclaim 

\demo{Proof}
Let $f_k = Q_k(f)$. Clearly, 

$$ \| \partial_x^{n/2 -1} f_k \|_{L_t^\infty L^2_x} \sim 
2^{k ( n/2 -1)} \| f_k \|_{L_t^\infty L^2_x}  \le \| f_k
\|_{\S^{(-1)}_k}.$$ Then (A.1)(ii) follows by taking $\ell^2$ norms both sides. 

To prove (A.1)(i) we proceed as follows. Given $q \ge 2$, let
$p^{\prime}$ be defined by $ \frac{1}{q} + \frac{n}{p'} =1$. 
In particular we have that $4 \le n \le p' \le 2\, n$. Since $n \ge 4$ we can now choose $ 2 \le r < p'$ such that $( q, r)$ is sharp admissible and 

$$ \| f_k \|_{L^q_t L_x^{p'}} \lesssim 2^{ k \gamma } \| f_k \|_{L^q_t
L^{r}_x}$$ by the Sobolev embedding where $\gamma$ is  given by 
and $\frac{1}{p'}= \frac{1}{ r } - \frac{\gamma}{n} $.

In particular then  $0 < \gamma =\frac{n}{r} -  \frac{n}{p'} =
\frac{n}{r}+ \frac{1}{q} - 1 $. From where by Lemma (2.7) we
can conclude that 
$$ \| f_k \|_{L^q_t L_x^{p'}} \lesssim \| f_k \|_{S_k^{(-1)}}. $$ 

The desired conclusion follows by taking $\ell^2$ norms both sides \qquad \qed . 

\enddemo
\vskip .1in

\proclaim{ (A.2) Lemma (First Multiplication estimate) } Let $ q  \ge 2$ and $p < 2$ such that $\dfrac{1}{q} + \dfrac{n}{p} = \dfrac{n}{2} +1 $. Then 
$$ P_{-1} : \S^{(-1)} \times \S^{(-1)} \lto L^q_t {\dot B}^{s}_{
{\tilde p}, 2} $$ for any $\tilde p \ge p $ and 
$ s =\frac{1}{q} + \frac{n}{\tilde p} - 1 .$

In particular note that when $q = 2$ we have that 
$$ \bigl( \sum_k 2^{2 k s} \| Q_k (\inv( f \cdot g)  \|_{L^2_t L^{\tilde
p}_x}^2 \bigr)^{1/2} \lesssim \| f\|_{\S^{(-1)} } \| g \|_{\S^{(-1)}}$$ 
 for any $ 2 \le \tilde p  <   2 n $ and $ s = \frac{n}{\tilde p} -
\frac{1}{2}$

\endproclaim

\demo{Proof}  
Let $f$ and $g$ be two function on $ \S^{(-1)}$. 

Let $q \ge 2$ and let $p^{\prime}$ and $p$ be defined by $$ \frac{1}{q} + \frac{n}{p'} =1 \quad \text{ and } \quad \frac{1}{p} = \frac{1}{2} + \frac{1}{p'}.$$

\noindent \underbar{Claim}: $$ \| P_{-1} ( f  \cdot g )
 \|_{ L^q_t \dot{W}^{n/2, p}} \lesssim  \| f \|_{ \S^{(-1)}}\| g \|_{ \S^{(-1)}}. $$

Assuming the claim we note that since $p < 2$ and $\frac{1}{q} + \frac{n}{p} = 1 + \frac{n}{2}$ 

$$ L^q_t \dot{W}^{n/2, p} \hookrightarrow   
L^q_t \dot{B}^{n/2}_{p, 2}  \hookrightarrow 
L^q_t \dot{B}^{s}_{{\tilde p}, 2} $$ 
provided  $s = \frac{1}{q} + \frac{n}{\tilde p}- 1$ and $\tilde p \ge
p $ as desired. 

To prove the claim we first note that by Lemma (A.1) we have, in
particular, the following two `endpoint estimates'.

$$ \|  f   \|_{L_t^\infty \dot{W}^{n/2-1, 2}} \lesssim  \| f
\|_{\S^{(-1)}}. \tag A.2)(i $$

$$  \| f   \|_{L^q_t L^{p'}_x } \lesssim 
 \| f   \|_{L^q_t \dot{B}^{0}_{p', 2}} \lesssim \| f
\|_{\S^{(-1)}}. \tag A.2)(ii $$ since $p' \ge 4 >2$. 

If $n=4$, the above two estimates suffice. For then,  
$$ \|  f \cdot g  \|_{L^q_t \dot{W}^{n/2 -1 , p }} \lesssim  \| f \|_{\S^{(-1)}} \| g \|_{\S^{(-1)}}\, \text{ where }\,  \frac{1}{p} = \frac{1}{2} + \frac{1}{p'} $$ since $p' \ge n$. In turn, this implies that 
$$ P_{-1} :   \S^{(-1)} \times \S^{(-1)} \lto  L^q_t \dot{W}^{n/2, p} $$
as desired. 

In the case $n \ge 5 $ however, we need to prove additional estimates. We consider separately the cases when 
$n$ is {\it even} first and then indicate the necessary modifications when $n$ is {\it odd}. 

More precisely, let $n \ge 5$ be even. Given $q, p'$ as above 
and $1\le j \le \dfrac{n}{2} - 2 $ let $ 0 < \theta_j < 1 $ be defined by 
$\theta_j = \dfrac{j}{(n/2 -1)}$. Next, let    $q_j $,  
$q_{n/2 -1 -j }$ and  $p'_j$,  $p'_{n/2 -1 -j }$ be any solutions to 
the following equations:
$$ \frac{1}{q_j} + \frac{n}{p'_j} =1 \qquad \text{ and } 
\qquad  \frac{1}{q_{n/2-1-j}} + \frac{n}{p'_{n/2-1-j}} =1  $$  

$$ \frac{1 - \theta_j}{q_j} + \frac{\theta_j}{ q_{n/2 -1 -j } } =
\frac{1}{q}\, , \qquad  \qquad  \frac{ 1-\theta_j}{p'_j} + \frac{\theta_j}{ p'_{n/2 -1 -j } } = \frac{1}{p'} $$

Next, let  $\dfrac{1}{{\tilde q}_j} =  \dfrac{1 - \theta_j}{q_j}$ and 
 $\dfrac{1}{{\tilde p'}_j} =  \dfrac{1 - \theta_j}{p'_j} + \dfrac{\theta_j}{2}$. 

\noindent We claim that 
 $$ \| \partial_x^{j} f \|_{L^{{\tilde q}_j}_t L^{{\tilde p'}_j }} \lesssim \| \partial_x^{j} f \|_{L^{{\tilde q}_j}_t \dot{B}^{0}_{{\tilde p'}_j , 2} }\lesssim 
\| f \|_{\S^{(-1)}}. \tag A.2)(iii $$

Indeed, the first inequality follows from the embeddings between Sobolev and Besov spaces since ${\tilde p'}_j \ge 2$. For the second one we have that $$ \align 
\| \partial_x^{j} f_k \|_{L^{{\tilde q}_j}_t L^{{\tilde p'}_j}}  &\sim
2^{k j} \| f_k \|_{L^{{\tilde q}_j}_t L^{{\tilde p'}_j}}\\ 
&= 2^{ k ( \theta_j (n/2 -1))}  \| f_k \|_{L^{{\tilde q}_j}_t L^{{\tilde p'}_j}} \\&=  2^{ k (\frac{1}{{\tilde q}_j} + \frac{n}{{\tilde p'}_j} -1 )}  \| f_k \|_{L^{{\tilde q}_j}_t L^{{\tilde p'}_j}} \lesssim \| f_k \|_{S_k^{(-1)}} \endalign $$    whence the second inequality follows by taking $\ell^2$-norms. 
   
\vskip .1in

Finally we put together with (A.2)(i), (A.2)(ii) and (A.2)(iii) to obtain that 
$$ \align  \|  f \cdot g  \|_{L^q_t \dot{W}^{n/2 -1 , p }} &= \| \partial_x^{n/2 -1} 
( f \cdot g ) \|_{L^q_t L_x^p} \\
&\lesssim \sum_{j=0}^{n/2-1} \| \partial_x^{j} f \|_{L^{{\tilde q}_{j}}_t L_x^{{\tilde p'}_{j}}}  \|\partial_x^{n/2 -1 -j} g \|_{L^{{\tilde q}_{n/2 -1 -j}}_t L_x^{{\tilde p'}_{n/2-1-j}}} \\
&\lesssim \| f \|_{\S^{(-1)}} \| g \|_{\S^{(-1)}}, \endalign  $$ as desired. 

\vskip .1in 

We indicate now the technicalities needed when $n \ge 5$ is odd. 
 Given $q, p'$ as above 
and $0\le j \le [\dfrac{n}{2}] - 1 $ let $ 0 < \theta_j < 1 $ to be
defined in a moment. As before, let  $q_j $,  
$q_{n/2 -1 -j }$ and  $p'_j$,  $p'_{n/2 -1 -j }$ 
be any solutions of the equations as above for $\theta_j$ and as
before let  $\dfrac{1}{{\tilde q}_j} =  \dfrac{1 - \theta_j}{q_j}$ and 
 $\dfrac{1}{{\tilde p'}_j} =  \dfrac{1 - \theta_j}{p'_j} +
\dfrac{\theta_j}{2}$. Now, 

$$ \align  \|  f \cdot g  \|_{L^q_t \dot{W}^{n/2 -1 , p }} &= \| \partial_x^{n/2 -1} 
( f \cdot g ) \|_{L^q_t L_x^p} \\
&\lesssim \sum_{j=0}^{[n/2]-1} \bigl( \| \partial_x^{j + 1/2}   f
\|_{L^{{\tilde q}_{j}}_t L_x^{{\tilde p'}_{j}}}  \|\partial_x^{n/2 -j
- 3/2} g \|_{L^{{\tilde q}_{n/2 -1 -j}}_t L_x^{{\tilde p'}_{n/2-1-j}}}
+ \\
& + \| \partial_x^{j }   f \|_{L^{{\tilde q}_{j}}_t L_x^{{\tilde
p'}_{j}}}\|\partial_x^{n/2 -j
- 1} g \|_{L^{{\tilde q}_{n/2 -1 -j}}_t L_x^{{\tilde p'}_{n/2-1-j}}} \bigr)
\endalign $$

\noindent For the first term inside the big sum we take $\theta_j =
\dfrac{j + 1/2}{(n/2 -1)}$; note that $\theta_{n/2-j-1} = 1 - \theta_j = \dfrac{n/2-j-3/2}{n/2 -1}$. Then for $0 \le j \le [n/2]-1 $ we have that
 $$ \| \partial_x^{j + 1/2}   f \|_{L^{{\tilde q}_j}_t L^{{\tilde p'}_j
}_x} \lesssim \| \partial_x^{j+ 1/2} f \|_{L^{{\tilde q}_j}_t
\dot{B}^{0}_{{\tilde p'}_j , 2} }\lesssim  \| f \|_{\S^{(-1)}}; $$
while 
$$ \| \partial_x^{n/2 -j -3/2}   g \|_{L^{{\tilde q}_{n/2-1-j}}_t L^{{\tilde p'}_{n/2-1-j}}_x} \lesssim \| \partial_x^{n/2-j-3/2} g \|_{L^{{\tilde q}_{n/2-j-1}}_t
\dot{B}^{0}_{{\tilde p'}_{n/2-j-1} , 2}} \lesssim  \| g \|_{\S^{(-1)}} $$

For the second we take $\theta_j = \dfrac{j}{(n/2 -1)}$ as in
(A.2)(iii), and the needed estimates follow just as in the even case. 

All in all we have that 
$$\|  f \cdot g  \|_{L^q_t \dot{W}^{n/2 -1 , p }}
\lesssim \| f \|_{\S^{(-1)}} \| g \|_{\S^{(-1)}},$$ as desired \qquad \qed.

\enddemo

\proclaim{ (A.3) Lemma (Second Multiplication estimate) }

$$ P_{-1} : \S^{(-1)} \times \S^{(-1)} \lto  \B_1 . $$

\endproclaim

\demo{Proof}  Let $f$ and $g$ be in $\S^{(-1)}$ and let $f_k =Q_k(f)$ and $g_j=Q_j(g)$ be their corresponding Littlewood-Paley projections. We write

$$ \align P_{-1} ( f \cdot g) &=   P_{-1} \bigl(  \sum_{ k, j \in \IZ} f_k \cdot g_j \bigr) \\ 
&=  P_{-1} \bigl(  \sum_{ k, j \in \IZ: k \ge j } f_k \cdot g_j \bigr) + P_{-1} \bigl(\sum_{ k, j \in \IZ : k < j} f_k \cdot g_j \bigr). \endalign  $$ By symmetry of the sums, it is enough to consider only one of them. The proof for the other is identical after exchanging $k$ and $j$. Hence we need to estimate 

$$\align & \sum_{l \in \IZ}  \| Q_l \bigl( P_{-1} \bigl(  \sum_{ k, j \in \IZ: k \ge j } f_k \cdot g_j \bigr)\bigr) \|_{L^1_t L^{\infty}_x} \\
&= \sum_{l \in \IZ} \| Q_l \bigl( P_{-1} \bigl(  \sum_{ k \in \IZ } \sum_{m \ge 0} f_k \cdot g_{k-m} \bigr) \bigr) \|_{L^1_t L^{\infty}_x} \endalign $$ 

Since $ \text{ supp } \widehat{(f_k \cdot g_{k-m})}\subseteq \{ \xi : |\xi| \le 2^k \} $ we have that $Q_l (f_k \cdot g_{k-m} ) \equiv 0 $ unless $k \ge l$. 

Therefore we can make the last sum less than or equal to 
$$\sum_{ l \in \IZ}  2^{-l} \sum_{ k \ge l } \sum_{m \ge 0} \| Q_l ( f_k \cdot g_{k-m} ) \|_{L^1_t L^{\infty}_x}. $$ On the other hand, we have that 
 $ \text{ supp } \widehat{(f_k \cdot g_{k-m})}  \cap \{ \xi : |\xi| << 2^{k-m} \} = \emptyset $ if $ m > 5$

Hence,  $Q_l (f_k \cdot g_{k-m} ) \equiv 0 $ unless $l=k$ and $m > 5$ or $ m \le 5$ and $ l < k $.  

We must then have that the above sum is  
$$ \lesssim \sum_{0 \le m \le 5} \sum_{ l \in \IZ}  2^{-l} \sum_{ k > l }\| Q_l ( f_k \cdot g_{k-m} ) \|_{L^1_t L^{\infty}_x} +  \sum_{m > 5} \sum_{ l \in \IZ}  2^{-l} \| Q_l ( f_l \cdot g_{l-m} ) \|_{L^1_t L^{\infty}_x}. $$
  
We consider the first sum first. 

$$ \align & \sum_{0 \le m \le 5} \sum_{ l \in \IZ}  2^{-l} \sum_{ k > l }\| Q_l ( f_k \cdot g_{k-m} ) \|_{L^1_t L^{\infty}_x}\\
& \lesssim  \sum_{0 \le m \le 5} \sum_{ l \in \IZ}  2^{-l} \sum_{ k > l } {(2^{n l})}^{\frac{n-3}{n-1}} \| f_k \|_{L^2_t L^{2p}_x} \| g_{k-m} \|_{L^2_t L^{2p}_x},\endalign $$  by Young's inequality with  $p = \frac{n-1}{n-3}$, $\frac{1}{p} + \frac{1}{p'} = 1 + \frac{1}{\infty} $  and  Cauchy-Schwartz inequality. 

\noindent The endpoint Strichartz estimates (2.1) now yield the bound 
$$\align & \sum_{0 \le m \le 5} \sum_{ l \in \IZ}  2^{-l} \sum_{ k > l } {(2^{n l})}^{\frac{n-3}{n-1}}  2^{ 2 k ( 1 + \frac{n}{(n-1)} - \frac{(n+1)}{2})} 
\| f_k\|_{\S^{(-1)}_k} \| g_{k-m}\|_{\S^{(-1)}_k} \\
& \sim \sum_{0\le m \le 5} \sum_{ k \in \IZ} \sum_{ l \le k } 2^{ l \frac{ w}{n-1}} 2^{- k  \frac{w}{n-1} } \| f_k\|_{\S^{(-1)}_k} \| g_{k-m}\|_{\S^{(-1)}_k} 
\endalign $$ 
where $w = (n-2)^2 - 3$ which is positive provided $n \ge 4$. Hence by summing first in $l$ and then applying Cauchy-Schwartz to the sum in $k$ we get that the above is $\lesssim \| f \|_{\S^{(-1)}}  \| g \|_{\S^{(-1)}}$ as desired. 

We proceed next with the second sum. 

$$ \align \sum_{m > 5} \sum_{ l \in \IZ}  2^{-l} \| Q_l ( f_l \cdot g_{l-m} ) \|_{L^1_t L^{\infty}_x} \\
&\lesssim  \sum_{m > 5} \sum_{ l \in \IZ}  2^{-l} 2^{n l/ r } \|  ( f_l \cdot g_{l-m} ) \|_{L^1_t L^{r}_x} \endalign $$ by Young's inequality with 
$r = 2 n $.   
Now, by H\"older's inequality we can bound the last sum by 

$$  \sum_{m > 5} \sum_{ l \in \IZ }  2^{-l} 2^{n l / r } \| f_l \|_{L^2_t L^{2r}_x} \|g_{l-m}  \|_{L^2_t L^{2r}_x}. $$
Since the pair $(2, 2r)$ is admissible we have by the Strichartz estimates that the above sum is up to a constant less than or equal to 
$$\align & \sum_{m > 5} \sum_{ l \in \IZ }  2^{-l} 2^{n l / r }
 2^{ 2 l ( 1/2 - n /(2 r) )} 2^{-m ( 1/2 - n/(2r))} \| f_l \|_{\S_l^{(-1)}} 
\|g_{l-m}  \|_{\S_l^{(-1)}} \\
&\lesssim \sum_{m > 5 } 2^{-m ( 1/2 - n/(2r))} \sum_{ l \in \IZ }  \| f_l \|_{\S_l^{(-1)}} \|g_{l-m}  \|_{\S_l^{(-1)}} \\
&\lesssim  \sum_{m > 5 } 2^{-m ( 1/2 - n/(2r))}\bigl( \sum_{ l \in \IZ } \,  \| f_l  \|^2_{\S_l^{(-1)}} \bigr)^{1/2} \bigl( \sum_{ l \in \IZ } \,  \|g_{l-m}  \|^2_{\S_l^{(-1)}} \bigr)^{1/2} \\
& \lesssim \| f \|_{ \S^{(-1)} }  \| g \|_{ \S^{(-1)} } 
\endalign $$ since by our choice of $r$,  $ 1/2 - n/(2r) = 1/4 > 0$ \qquad \qed.  

\enddemo

\newpage
\enddocument